# ON THE PATH DENSITY OF A GRADIENT FIELD

By Christopher R. Genovese, Marco Perone-Pacifico,
Isabella Verdinelli and Larry Wasserman

*Carnegie Mellon University, University of Rome, University of Rome and Carnegie Mellon University*

We consider the problem of reliably finding filaments in point clouds. Realistic data sets often have numerous filaments of various sizes and shapes. Statistical techniques exist for finding one (or a few) filaments but these methods do not handle noisy data sets with many filaments. Other methods can be found in the astronomy literature but they do not have rigorous statistical guarantees. We propose the following method. Starting at each data point we construct the steepest ascent path along a kernel density estimator. We locate filaments by finding regions where these paths are highly concentrated. Formally, we define the density of these paths and we construct a consistent estimator of this path density.

**1. Introduction.** The motivation for this paper stems from the problem of finding filaments from point process data. Filaments are one-dimensional curves embedded in a point process or random field. Identifying filamentary structures is an important problem in many applications. In medical imaging, filaments arise as networks of blood vessels in tissue and need to be identified and mapped. In remote sensing, river systems and road networks are common filamentary structures of critical importance [Lacoste, Descombes and Zerubia (2005) and Stoica, Descombes and Zerubia (2004)]. In seismology, the concentration of earthquake epicenters traces the filamentary network of fault lines. This paper is motivated by a cosmological application: the detection of filaments of matter in the universe.

Figures 1 and 2 show two examples. The first example is a cosmology data set showing positions of galaxies. The second is a synthetic example. In each case, the upper left plot shows the data, which exhibits an apparent filamentary structure. A density estimate reveals this structure more clearly.

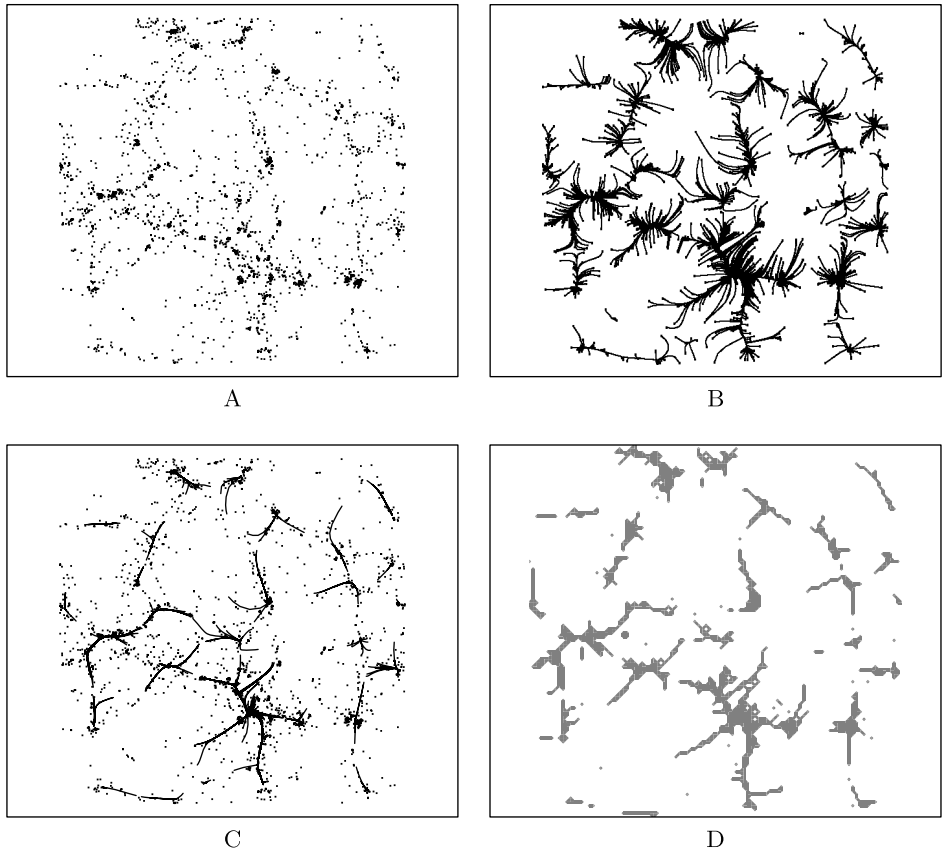

FIG. 1. *Cosmological data: data set (plot* A*), steepest ascent paths of all data points (plot* B*), paths after cutting the first five iterations (plot* C*) and level set at* 90% *quantile of the estimated path density (plot* D*).*

In particular, the steepest ascent paths of the density estimate—the paths from each point that follows the gradient—tend to concentrate along the filaments. The upper right plot in each figure shows the collection of paths, and the lower left plot trims off the early iterations. The result, on the bottom right in the figures, shows a high density of paths around the filaments. The empirical observation that the steepest ascent paths concentrate around the filaments motivates our approach in this paper. Specifically, we characterize the density of steepest ascent path—what we call the *path density* below—and construct a consistent estimator of the path density using the steepest ascent paths of a density estimator. Our purpose in this paper is to explore this idea in detail. We apply this idea to the filament problem by showing that the path density is high near filaments and that the set where the estimated path density is large essentially capture the filaments.

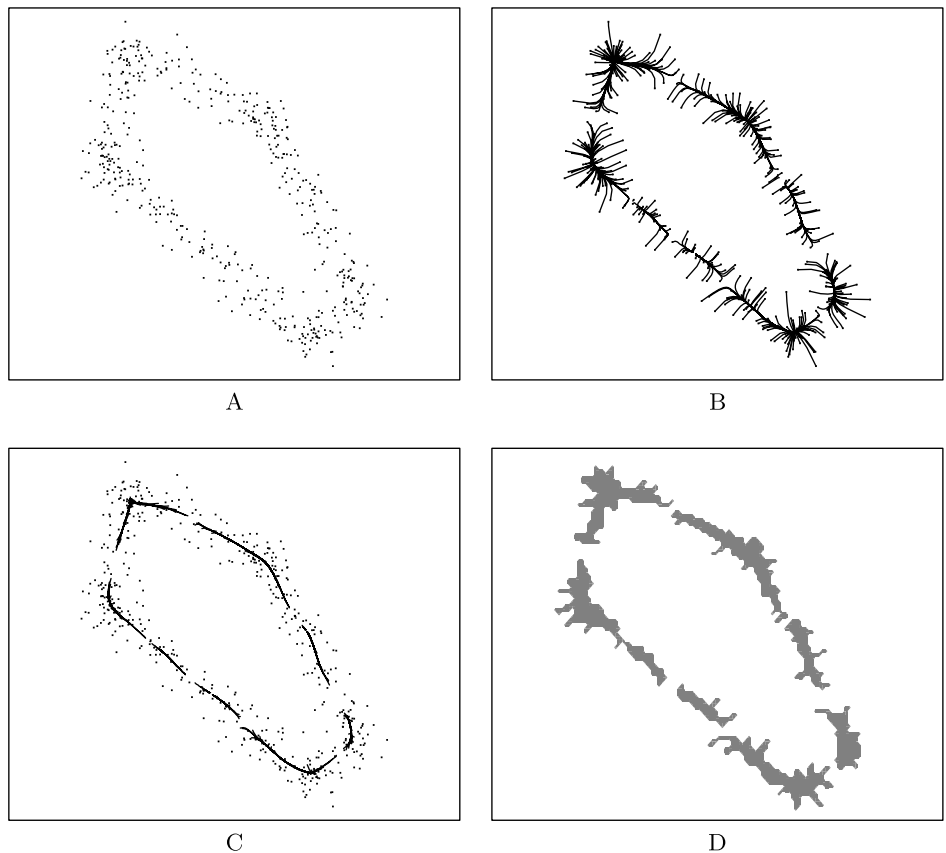

Fig. 2. *Simulated example: data set (plot* A*), steepest ascent paths of all data points (plot* B*), paths after cutting the first three iterations (plot* C*) and level set at* 90% *quantile of the estimated path density (plot* D*).*

There are existing techniques for finding one (or a few) filaments. Examples from statistics and machine learning include Arias-Castro, Donoho and Huo (2006), Kegl et al. (2000), Sandilya and Kulkarni (2002) and Tibshirani (1992). But none of these methods are practical for noisy data sets with large numbers of complicated filaments. Other methods can be found in the astronomy literature: see Stoica et al. (2005), Eriksen et al. (2004), Novikov, Colombi and Doré (2006), Sousbie et al. (2006) and Barrow, Bhavsar and Sonoda (1985). But these methods do not have rigorous statistical guarantees. Thus the problem of reliably finding many filaments simultaneously remains largely unsolved. Hence our search for new methods. We would like to emphasize that the methods proposed here are still very far from being a complete methodology for filament analysis. Rather our results are a modest first step toward a novel and promising methdology. A recent paper on estimating integral curves is Koltchinskii (2007), although the setting and methodology are quite different.

Let us now give a heuristic description of the basic problem we study. Let $X_1, \ldots, X_n$ be a sample from a distribution $\mu_X$ with density $g_X$. We estimate the function

$$p(x) = \lim_{r \to 0} \frac{\mathbb{P}(P(X) \text{ intersects } B(x,r))}{r},$$

where $P(X)$ is the steepest ascent path defined by $g_X$ starting at $X$ and $B(x,r) = \{y : \|y - x\| < r\}$. We use balls here because they are simple, but any suitably rich collection of open neighborhoods yields the same function.

Our main results are as follows. Theorem 6 explicitly characterizes the path density $p(x)$. Theorem 3 constructs an estimator $\widehat{p}_n$ such that

$$\sup_x |\widehat{p}_n(x) - p(x)| = O_P\left(\sqrt{\frac{\log(1/h_n)}{nh_n^4}}\right) + O_P\left(\sqrt{\frac{\log n}{n\nu_n}}\right) + O(\nu_n) + O(h_n^2),$$

where $h_n$ and $\nu_n$ are appropriately chosen bandwidth parameters. We then apply these ideas to the problem of identifying filaments from point-process data and show, in Theorem 4, that the level sets of the path density concentrate around the filaments.

1.1. *Notation.* We assume a probability space $(\Omega, \mathcal{S}, \mathbb{P})$ on which all our random variables are defined. We consider random variables taking values in $\mathbb{R}^2$. Throughout the paper, $U_0$ denotes a large, fixed compact subset of $\mathbb{R}^2$ whose properties are specified later.

Let $B(x,r)$ denote the open ball of radius $r$ centered on $x$. Denote the closure of $B(x,r)$ by $\overline{B}(x,r)$ and its boundary by $\partial B(x,r)$. For any set $A$, define the $r$-dilation by

$$B(A,r) = \bigcup_{x \in A} B(x,r). \tag{1}$$

We use $1_A(x)$ as the indicator function of the set $A$.

If $f : \mathbb{R}^n \to \mathbb{R}^m$, then we use $Df$, $D^2 f$, $D_i f$ to denote various total and partial derivatives. Specifically, $Df$ is a linear map from $\mathbb{R}^n \to \mathbb{R}^m$, $D^2 f$ denotes the Hessian matrix of $f$, and $D_i f$ is the partial derivative of $f$ with respect to the $i$th argument. When $m = 1$, $Df$ is a gradient, which is convenient to view as a vector in $\mathbb{R}^n$. For this purpose, we use $\nabla f$ to denote the column vector $(Df)^T$.

We use $\phi$ to denote the univariate standard normal density and $\Phi$ the corresponding distribution function. For $\sigma > 0$, $\phi_\sigma$ and $\Phi_\sigma$ are the corresponding $N(0, \sigma^2)$ functions. In $\mathbb{R}^d$ for $d \geq 2$, we use $\varphi$ for the $d$-dimensional standard normal density and $\varphi_\sigma$ for the $N(0, \text{diag}(\sigma^2, \ldots, \sigma^2))$ density.

1.2. *Outline.* In Section 2, we define and characterize the path density function. In Section 3, we define an estimator of the path density function and find its rate of convergence. Section 4 describes the problem of identifying filaments from point-process data, showing that the path density is large near filaments and small elsewhere. Proofs are relegated to Sections 5 through 7. We close with some general remarks in Section 8.

## 2. Integral curves and path densities.

2.1. *Flows and the local group.* If $V$ is a smooth ($C^\zeta$ for $\zeta \geq 1$) vector field on $\mathbb{R}^2$, then we might imagine putting a test particle at any $x \in \mathbb{R}^2$ and letting it flow with velocity given by the vector field.

It is known [the Picard–Lindelof theorem, Irwin (2001)] that this idea is well defined. The paths followed by the test particle starting at different points fit together in a consistent way. In particular, in any neighborhood $U$ of $x \in \mathbb{R}^2$, there is a neighborhood $U_1 \subset U$, an interval $I \subset \mathbb{R}$ containing 0 and a $C^\zeta$ mapping $\psi: I \times U_1 \to U$ such that (i) $\psi(0, x) = x$, (ii) $\frac{\partial}{\partial t}\psi(t, x) = V(\psi(t, x))$, and (iii) if $s, t, s + t \in I$, $\psi(s + t, x) = \psi(s, \psi(t, x))$. For each $t$, we define a mapping $\psi^t x \equiv \psi(t, x)$ that gives the point obtained by following the flow from $x$ for time $t$. In these terms, property (iii) becomes $\psi^{s+t} = \psi^s \circ \psi^t$, giving a group-like structure with composition as the product. Because of this, the mappings $\psi^t$ are called the local one parameter group of diffeomorphisms generated by $V$. Moreover, the paths $t \mapsto \psi^t x$, called *integral curves* or local flows of the vector field, are unique in the sense that integral curves are equal where their domains overlap. Thus, the integral curves create an equivalence relation on $\mathbb{R}^2$, where two points are equivalent if they are on the same integral curve.

In certain cases, these local flows can be extended to a *global* flow, a mapping $\psi: \mathbb{R} \times \mathbb{R}^2 \to \mathbb{R}^2$ satisfying properties (i), (ii) and (iii) above with $I = \mathbb{R}$. The mappings $\psi^t: \mathbb{R}^2 \to \mathbb{R}^2$ are $C^\zeta$ diffeomorphisms that form a group under composition. Such global flows exist, for instance, when the domain of the vector field is compact, when $V$ has compact support or when the domain is a Banach space [see Theorem 3.39 in Irwin (2001)], which is the case for $\mathbb{R}^2$.

We derive a vector field from the gradient of the density $g_X$. A critical point of $g_X$ is one where the gradient of $g_X$ equals 0. Any other point is called a regular point of $g_X$. We make the following assumption in what follows.

ASSUMPTION 1. $g_X$ is a $C^{\zeta+1}$ function on $\mathbb{R}^2$ for $2 < \zeta \leq \infty$.

ASSUMPTION 2. All the critical points of $g_X$ are nondegenerate, meaning that the Hessian is nonsingular.

ASSUMPTION 3.

$$\lim_{d \to \infty} \inf_{\substack{B \\ \mathrm{diam}(B)=d}} \sup_{\partial B} \left| \frac{\nabla g_X}{\|\nabla g_X\|} - \mathbf{n}_B \right| = 0, \tag{2}$$

where the infimum is over closed balls in $\mathbb{R}^2$ of given diameter $d$ and where $\mathbf{n}_B$ is the outward pointing, unit normal vector field defined on $\partial B$.

The importance of Assumption 2 is that the behavior of $g_X$ around nondegenerate critical points is locally quadratic. See Remark 3 below. Nondegenerate critical points are necessarily isolated. Moreover, by Assumption 3, all $g_X$'s critical points lie in a compact set. These facts together with Assumption 2 imply that the critical points of $g_X$ are finite in number.

Assumption 3 ensures that the gradient of $g_X$ is nearly radial far enough out from some point. This condition is satisfied by a wide variety of common distributions whose critical points lie within a compact set. For example, in Section 4, we show that the results extend to quite general mixtures of normal distributions. A stronger alternative assumption that is easier to understand is that $g_X$ has compact support, which is certainly sufficient for the results that follow.

From the assumptions above it follows that $V = \nabla g_X$ gives a $C^\zeta$-vector field on $\mathbb{R}^2$ and therefore generates a unique global flow $\psi$. Here, the flow $\psi^t x$ moves along the direction of steepest ascent. If the support of $g_X$ is not necessarily compact, as in Section 4, we can restrict to a compact set $U_0$ containing all critical points of $g_X$. Restriction to $U_0$ requires that we define an interval $I_x = [a_x, b_x]$ such that $\psi^t x \in U_0$ whenever $t \in I_x$. If Assumption 3 holds, the interval is the entire real line, but we maintain the interval notation in the proofs to facilitate extension of the results in later sections.

The global flow $\psi$ on $\mathbb{R}^2$ generates an equivalence relation on $U_0$ (or $\mathbb{R}^2$ for that matter). For $x, y \in U_0$, say that $x \sim y$ if $\psi^t y = x$ for some $t \in \mathbb{R}$. We say that $x$ precedes $y$, $x \preceq y$ ($x \prec y$), if $\psi^{-t} y = x$ for $t \geq 0$ ($t > 0$) and $y$ succeeds $x$, $y \succeq x$ ($y \succ x$), if $\psi^t y = x$ for $t \geq 0$ ($t > 0$). Here precedence and succession refer to the flow in directions of increase of $g_X$, that is, along the vector field $V$. For any $A \subset U_0$, define the *reverse evolution* of $A$ under the flow by

$$\mathcal{V}(A) = \{y \in U_0 : y \preceq x, x \in A\}, \tag{3}$$

the set of points in $U_0$ that precede a point in $A$. With minor abuse of notation, we also write $\mathcal{V}(x) = \mathcal{V}(\{x\})$.

A simple example may clarify these definitions. Let $g(x) = -\frac{1}{2}\|x\|^2$, a simple quadratic. The gradient of $g$ is a vector field that at each point $x$ gives the direction of steepest ascent $\nabla g(x) = -x$. This vector field specifies differential equations at every point $x : \dot{\gamma}(t) = -\gamma(t)$ with $\gamma(0) = x$, where $\gamma = (\gamma_1, \gamma_2)$ is a curve, that is, a function from $\mathbb{R}$ to $\mathbb{R}^2$. This is the integral curve. The mapping $\psi$ packages together all these integral curves, and in this example, $\psi^t x = (x_1 e^{-t}, x_2 e^{-t})$. Figure 3 illustrates the flows in a generic case.

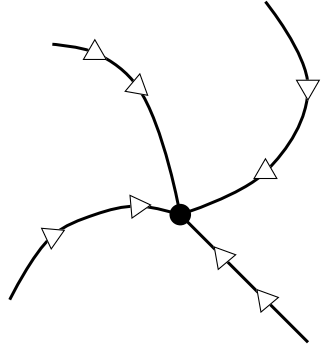

Fig. 3. *Flow lines in a typical gradient field, moving toward a local maximum.*

As we will see in Theorem 4, the integral curves of the flow $\psi$ concentrate near the filaments. We would thus expect that following the flow induced by the gradient of $g_X$ will carry a collection of random points near to the filaments. This is the idea behind the estimators described in the first section. We quantify this concentration through the *path measure* $\pi$, which we define to be the probability that the flow from a random point hits a given set. That is, for $A \subset U_0$,

$$\pi(A) = \mu_X(\mathcal{V}(A)). \tag{4}$$

This is a sub-additive set function (in fact, an infinitely-alternating Choquet capacity), and does not have a density in the Radon–Nikodym sense. However, we show below that the following density is well defined.

DEFINITION 1. The path density associated with $g_X$ is the function $p \colon \mathbb{R}^2 \to [0, \infty]$ defined by

$$p(x) = \lim_{r \to 0} \frac{\pi(B(x,r))}{r}. \tag{5}$$

REMARK 1. One might have expected the denominator in the expression above to be $r^2$ since we are working in $\mathbb{R}^2$ but, as it turns out, $\pi(B(x,r))$ is order $r$.

REMARK 2. To clarify, $\pi$ is a set function with $\pi(A)$ giving the probability that for a random point drawn from $g_X$, the flow from that point hits $A$. The path density $p$ is a function on $\mathbb{R}^2$ with $p(x)$ describing the $\pi$-probability in infinitesimal neighborhoods of $x$. See Figure 4.

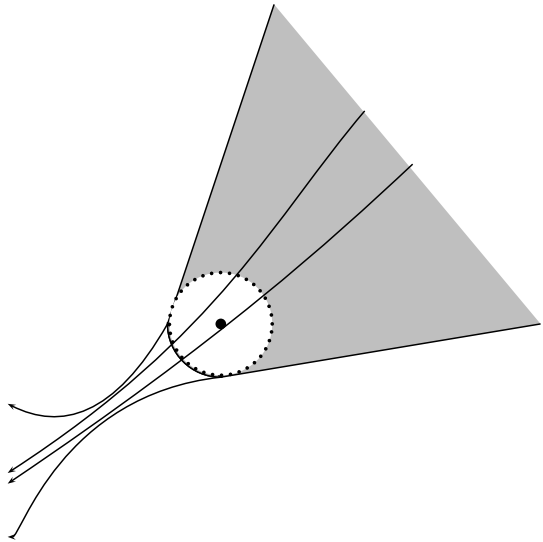

Fig. 4. *A set A represented by the white disk and the set of points whose flow lines hit A, the union of A and the shaded gray region. $\pi(A)$ is the probability content of the latter set. The path density is the limit of these probabilities, suitably normalized, as A shrinks.*

The next theorem gives the main features of the path density. A more detailed characterization of $p$ is given by Theorem 6 in Section 6. Let $\mathcal{M}$ denote the set of local maxima of $g_X$ and $\mathcal{H}$ denote the set of saddlepoints of $g_X$.

THEOREM 1. *Under Assumptions 1–3, the path density $p$ is an upper-semicontinuous function with the following properties:*

1. *if $x$ is a local minimum of $g_X$, $p(x) = 0$;*
2. *if $x$ is a local maximum of $g_X$, $p(x) = \infty$;*
3. *$p$ is continuous on $(\mathcal{M} \cup \mathcal{H})^c$ and bounded on $\mathcal{M}^c$.*

REMARK 3. The result of Theorem 1 depends on the nondegeneracy of $g_X$'s critical points. While nondegenerate critical points are isolated, degenerate critical points need not be. For example, a density with a linear ridge can produce a line of local maxima. A mixture of normals over two closely spaced, parallel lines would generate a similar ridge between the lines. In such cases, the local behavior of the function around critical points need not be local quadratic as it must around nondegenerate critical points. Near a ridge, for instance, the flow lines from two neighboring points will proceed in parallel toward the ridge rather than converging on a local maximum. As such, $p$ is finite at any point on the ridge, although $\pi(B(R,r))/r \to \infty$ as $r \to 0$ where $R$ is the set of points on the ridge.

The good news is that degeneracy is rare in the sense that functions without degenerate critical points are dense in the space of smooth functions. The proof of this follows from Sard's theorem by adding to $g_X$ a linear function that is arbitrarily small over $U_0$; the linear function changes the locations of the critical points without changing the Hessian and almost every linear function eliminates the nondegeneracy. [See Guillemin and Pollack (1974) for further discussion of this point.]

The next theorem gives a more precise approximation of the path measure in terms of the path density. It will be used extensively in what follows.

THEOREM 2. *Let $x$ be a regular point of $g_X$. Then, the following hold:*

1. *As $r \to 0$,*

(6) $$\pi(B(x,r)) = rp(x) + O(r^2).$$

2. *There exists a curve $\gamma$, parameterized by arc length, such that*

(7) $$\pi(B(x,r)) = \frac{1}{2}\int_{-r}^{r} p(\gamma(t))\,dt + O(r^2).$$

**3. Estimating the path density.** Recall that $X_1, \ldots, X_n$ are independent observations from the distribution $\mu_X$ with density $g_X$. Our goal is to estimate the path density $p$.

First we estimate $g_X$ using the kernel estimator

$$\widehat{g}_n(x) = \frac{1}{n}\sum_{i=1}^{n} \frac{1}{h^2} K\left(\frac{\|x - X_i\|}{h}\right).$$

Given an arbitrary point $x \in \mathbb{R}^2$, let $P(x) = \{\psi^t x : 0 \leq t \leq \infty\}$ denote the integral curve starting at $x$ and let $\widehat{P}(x)$ denote the integral curve obtained by replacing $g_X$ with $\widehat{g}_n$. [Recall that $\psi^t x$ is the point obtained by starting at $x$ at time 0 and following the vector field until time $t$. The integral curve $P(x)$ shows the whole path going forward in time.] We define the *path density estimator*

(8) $$\widehat{p}_n(x) = \frac{1}{n}\sum_{i=1}^{n} \frac{1}{\nu} K\left(\frac{\inf_{z \in \widehat{P}(X_i)} \|x - z\|}{\nu}\right).$$

Thus, $\widehat{p}_n(x)$ is essentially a weighted count of how many observed paths $\widehat{P}(X_1), \ldots, \widehat{P}(X_n)$ get close to $x$. It is not necessary to use the same kernel $K$ for $\widehat{g}_n$ and $\widehat{p}_n$, but for simplicity we shall take them to be the same up to a normalizing constant.

REMARK 4. A numerical approximation to $\widehat{P}(x)$ can be obtained using the mean shift algorithm [Fukunaga and Hostetler (1975) and Cheng (1995)]. Given an arbitrary point $x \in \mathbb{R}^2$, the mean shift algorithm defines a sequence $\widehat{s}(x) = (x^0 = x, x^1, x^2, \ldots)$, where

$$(9) \qquad x^{k+1} = \frac{\sum_{i=1}^n X_i K(1/h \| x^k - X_i \|)}{\sum_{i=1}^n K(1/h \| x^k - X_i \|)}.$$

The sequence $\widehat{s}(x)$ converges to a mode of $\widehat{g}_n$. Conversely, for each mode of $\widehat{g}_n$ there exists a point $x$ such that $\widehat{s}(x)$ converges to that mode. A smooth interpolation of the points in $\widehat{s}(x)$ can be regarded as a numerical approximation to $\widehat{P}(x)$.

We make the following assumptions about the kernel $K$.

(K1) $K : [0, \infty) \to [0, \infty)$ is a nonincreasing, bounded, square integrable function and has bounded derivative.

(K2) $K$ satisfies the Giné and Guillou (2002) conditions, namely, $K$ belongs to the linear span of functions with the following property: the set $\{(s, u) : K(s) \geq u\}$ can be represented as a finite number of Boolean operations among sets of the form $\{(s, u) : a(s, u) \geq b(u)\}$ where $a$ is a polynomial and $b$ is any real function.

(K3) $K$ satisfies the following tail condition: as $x \to \infty$

$$K(x) = O(xe^{-x}).$$

Condition (K2) is somewhat abstract but, unfortunately, such a technical condition appears to be necessary. The role (K2) plays in the proof is to ensures that the VC-condition holds, thus allowing uniform control of an empirical process. The good news, as Giné and Guillou (2002) point out, is that (K2) is satisfied by most common kernels.

THEOREM 3. *Suppose that $K$ satisfies* (K1)–(K3), *and that the bandwidths $h_n$ and $\nu_n$ satisfy the following conditions:*

$$h_n \to 0, \qquad \frac{nh_n^4}{|\log h_n|} \to \infty, \qquad \frac{|\log h_n|}{\log \log n} \to \infty, \qquad h_n^2 \leq c_1 h_{2n}^2,$$

*for some $c_1$, and*

$$\nu_n \to 0, \qquad \frac{n\nu_n}{|\log \nu_n|} \to \infty, \qquad \frac{|\log \nu_n|}{\log \log n} \to \infty, \qquad \nu_n \leq c_2 \nu_{2n},$$

*for some $c_2$. Further, assume that $g_X$ is bounded. Then,*

$$(10) \qquad \sup_x |\widehat{p}_n(x) - p(x)| = O_P\left(\sqrt{\frac{\log(1/h_n)}{nh_n^4}}\right) + O_P\left(\sqrt{\frac{\log n}{n\nu_n}}\right) + O(\nu_n) + O(h_n^2),$$

*where the $\sup_x$ is taken over all regular points.*

The best rate is obtained by setting $h_n \asymp (\log n)^{1/4}/n^{1/8}$ and $\nu_n \asymp \log n/n^{1/3}$.

COROLLARY 1. *Setting $h_n \asymp (\log n)^{1/4}/n^{1/8}$ and $\nu_n \asymp \log n/n^{1/3}$ we have*

$$\sup_x |\widehat{p}_n(x) - p(x)| = O_P\bigg(\frac{\sqrt{\log n}}{n^{1/4}}\bigg), \tag{11}$$

*where the $\sup_x$ is taken over all regular points.*

**4. Filament detection.** Galaxies, being large and bright and having a tendency to cluster together, serve as tracers of matter. At large enough scales, the universe looks the same in every direction, so astronomers were surprised when their maps of the galaxy locations revealed complicated, systematic structures—clusters, walls, sheets, and voids. But the most striking feature in these maps is the vast network of filaments, often called the "cosmic web." Panel A in Figure 1 shows a map from one such survey, with the galaxy locations selected from a two-dimensional slice of the universe. Astronomers want to identify the filaments to help them characterize this large scale structure and in turn constrain the physics of the universe's evolution.

Astronomers have substantial literature on the problem of estimating filaments. Early efforts involved the use of standard spatial techniques, including high-order correlations, shape statistics from Luo and Vishniac (1995) and Minkowski functionals. [See the book by Martinez and Saar (2002) for an overview of such approaches.] Barrow, Bhavsar and Sonoda (1985) applied minimal spanning trees to the problem, but these trees give primarily local descriptors of neighorhood relations. The skeleton method [Stoica et al. (2005), Eriksen et al. (2004), Novikov, Colombi and Doré (2006) and Sousbie et al. (2006)] involves smoothing the point distribution to estimate the skeleton: a network connecting saddlepoints to local maxima of the point density with edges parallel to the gradient of the density. Stoica et al. (2005) develops an automated method for tracing a filamentary network based on marked point processes; Stoica, Martinez and Saar (2007) extends this to three-dimensional networks.

In the statistics literature, filaments are similar to, but distinct from, principal curves [Hastie and Stuetzle (1989)]. Hastie and Stuetzle (1989) define a principal curve $c$ by the self-consistency relation $c(t) = \mathbb{E}(X|\pi_c(X) = t)$ where $\pi_c(x)$ is the projection of $x$ onto $c$. Modified definitions are given by Kegl et al. (2000) and Sandilya and Kulkarni (2002). They call $c$ a principal curve if $c$ minimizes $\mathbb{E}\|X - \pi_c(X)\|^2$ subject to $c$ lying in a prespecified set of curves such as all curves of length bounded by some constant or with bounded curvature. In either case, they show that there exist estimators of a filament with convergence rate $O(n^{-1/3})$. However, as noted by Hastie and Stuetzle (1989), the principal curve $c$ is not equal to the filament.

Recently, Arias-Castro, Donoho and Huo (2006), considered the problem of finding a single filament in a noisy background. They obtain precise minimax bounds on how sparse a filament can be and still be detectable. However, their method does not easily extend to the case where there are many filaments. Also, they can only detect the presence of a filament but they do not produce an estimate of the filament itself.

4.1. *A model for filaments.* We assume that the data are a realization of an inhomogeneous Poisson process on $U_0$ such that conditional on the number of observations, $n$, we have $n$ i.i.d. draws $X_1, \ldots, X_n$ from the density $g_X$ on $U_0$. We model $g_X$ as a mixture of three types of components: filaments, clusters and the background.

A filament is a smooth, nonself-intersecting curve $f:[a,b] \to U_0$. We consider a finite collection of $C^2$ filaments $f_1, \ldots, f_{m_F}$ of lengths $\ell_1, \ldots, \ell_{m_F}$, which are allowed to be positioned arbitrarily within $U_0$, including the possibility that distinct filaments intersect. Because our inferences are independent of how the curve is parameterized, we parameterize the curve by arclength, taking $a = 0$ and $b = \ell(f)$, the length of $f$ in $\mathbb{R}^2$. Define $\mathcal{F} = \bigcup_{i=1}^{m_F} \mathrm{range}(f_i)$ to be the filament set. Clusters are concentrations of points around a cluster center $z$. We consider a finite collection of clusters which we denote $\mathcal{C} = \{z_{m_F+1}, \ldots, z_m\}$. The background generates points that are not on filaments or in clusters, and we model it as a homogeneous Poisson process.

Specifically, $g_X$ takes the following form:

$$
g_X(x) = \alpha_0 \frac{1_{U_0}(x)}{\mathcal{L}(U_0)} + \sum_{i=1}^{m_F} \alpha_i \int_0^{\ell_i} w_i(s) \varphi_{\sigma_i}(x - f_i(s))\, ds \\
+ \sum_{i=m_F+1}^{m} \alpha_i \varphi_{\sigma_i}(x - z_i),
\tag{12}
$$

where $\mathcal{L}$ denotes the Lebesgue measure and each $w_i$, $i = 1, \ldots, m_F$, is a strictly positive probability density on $[0, \ell_i]$; $\sigma_1, \ldots, \sigma_m > 0$, and $0 \leq \alpha_1, \ldots, \alpha_m \leq 1$ with $\sum_i \alpha_i = 1$. Recall that $\varphi_\sigma$ denotes the symmetric normal distribution over $\mathbb{R}^2$. We define $\mu_X$ to be the probability measure on $\mathbb{R}^2$ with density $g_X$. Generally we take all $\sigma_i$ equal to a common value $\sigma$.

Our model for the filament data is a generalization of the model in Tibshirani (1992): select a random point along the filament and generate from a normal centered on that point. Note that clusters are zero-length filaments, with the corresponding $w_i$ equal to a delta function at 0.

When $\alpha_0 = 0$, we say that $g_X$ is *background free*. In this case, $g_X$ is a $C^\infty$ function on all of $\mathbb{R}^2$. Our results from Section 2 apply to the background-free case, but they extend directly to the case where $\alpha_0 > 0$ because $\nabla g_X$ is a scaled version of the background-free gradient on $U_0$. We can either restrict the global flow to $U_0$ or adapt it to $U_0$ by constructing a $C^\infty$ function $\omega$ with compact support that equals 1 on $U_0$ and is zero outside a small dilation of $U_0$. Our results carry over immediately when $\mathbb{R}^2$ is replaced by $S^2$, the two sphere, because $S^2$ is a compact, two-dimensional manifold.

Under model (12), we can without loss of generality take $U_0$ large enough to contain all the critical points. This is true because, in the background-free case, solving $\nabla g_X(x) = 0$ gives

$$
\begin{aligned}
x = &\left( \sum_{i=1}^{m_F} \alpha_i \int_0^{\ell_i} f_i(s) w_i(s) \varphi_{\sigma_i}(x - f_i(s))\, ds + \sum_{i=m_F+1}^{m} z_i \alpha_i \varphi(\sigma_i(x - z_i)) \right) \\
&\times \left( \sum_{i=1}^{m_F} \alpha_i \int_0^{\ell_i} w_i(s) \varphi_{\sigma_i}(x - f_i(s))\, ds + \sum_{i=m_F+1}^{m} \alpha_i \varphi(\sigma_i(x - z_i)) \right)^{-1},
\end{aligned}
\tag{13}
$$

where the terms come from (12), and so $x$ must lie in the convex hull of $\mathcal{F} \cup \mathcal{C}$, which is compact. (A background does not change the gradient on $U_0$, leading to the same conclusion.) By Assumption 2, all critical points of $g_X$ are isolated, and therefore $\mathcal{M}$ and $\mathcal{H}$ and the set of local minima are all finite.

Note that under our filament model, $g_X$ does not have compact support. However, because the filaments are contained in a compact set, $g_X$ has normal tails in all directions outside large enough balls. In particular, $\|\frac{\nabla g_X}{\|\nabla g_X\|} - \mathbf{n}\| \to 0$ uniformly as the ball radius increases, where $\mathbf{n}$ is the outward pointing normal to the ball's boundary. This ensures that any contribution to the path density from points outside $U_0$ is negligible.

4.2. *Capturing filaments.* Now we show that the flow lines concentrate near the filaments. We do that by showing that the sets on which $p$ exceeds some value $\lambda$ lies in a set that is close to the filaments and getting closer as $\lambda$ grows. Two complications are the local maxima, at which $p(x) = \infty$ and thus exceeds any $\lambda$, and saddle points at which $p(x)$ is roughly four times the value of $p$ in a neighborhood around $x$. Define $\sigma = \max_i \sigma_i$ and

$$
d(\lambda) = \sigma \sqrt{2 \log(1/2\pi\sigma^2 \lambda)}.
\tag{14}
$$

THEOREM 4. *Let $\mathcal{A} = \mathcal{F} \cup \mathcal{C} \subset U_0$, the set of points that are on a filament or at a cluster center. Let $\mathcal{H}$ denote the set of saddle points of $g_X$. Define*

$$
\mathcal{E}_\lambda = \{x \in U_0 : p(x) > \lambda\}.
\tag{15}
$$

Then, for all $\varepsilon > 0$ and all $\lambda \geq \varepsilon$,

(16) $\quad\quad\quad \mathcal{E}_\lambda \cap \mathcal{M}^c \subset B(\mathcal{A}, d(\lambda) + \varepsilon) \cup (\mathcal{H} \cap B(\mathcal{A}, d(4\lambda) + \varepsilon)),$

(17) $\quad\quad \mathcal{E}_\lambda \cap \mathcal{M}^c \cap \mathcal{H}^c \subset B(\mathcal{A}, d(\lambda) + \varepsilon).$

Theorems 1 and 4 together imply that $\{x : \lambda < p(x) < \infty\} \subset B(\mathcal{A}, d(4\lambda) + \varepsilon)$ and that $\{x : \lambda < p(x) < \infty\} \cap \mathcal{H}^c \subset B(\mathcal{A}, d(\lambda) + \varepsilon).$

COROLLARY 2. *There exists $\mathcal{A}_\sigma$ such that $\mathcal{A}_\sigma \subset \mathcal{A}$,*

$$\mathcal{A}_\sigma \subset \mathcal{E}_\lambda \subset B(\mathcal{A}_\sigma, d(\lambda))$$

*and $d_H(\mathcal{A}_\sigma, \mathcal{A}) \to 0$ as $\sigma \to 0$. Hence, for every $\xi > 0$ there exists $\lambda$ such that $d_H(\mathcal{A}, \mathcal{E}_\lambda) < \xi + O(\sigma)$.*

This corollary says that, as long as the noise level $\sigma$ is small, the level sets approximate the filaments.

We make use of results in Cuevas and Fraiman (1997). A set $S$ is *standard* if for every $\lambda > 0$ there exists $0 < \delta < 1$ such that

$$\mathcal{L}(B(x, \varepsilon) \cap S) \geq \delta \mathcal{L}(B(x, \varepsilon))$$

for every $x \in S$, $0 < \varepsilon \leq \lambda$.

THEOREM 5. *Suppose that $K$ satisfies conditions* (K1)–(K3), *and that $\mathcal{E}_\lambda$ is standard. Let $\widehat{\mathcal{F}} = \{x : \widehat{p}(x) \geq \lambda + c_n\}$ where $c_n \geq 0$ and $c_n \to 0$. Suppose that $\beta_n \to \infty$, $\beta_n \nu_n \to \infty$, and $\beta_n \nu_n / c_n$ is bounded. Then,*

$$\beta_n d(\mathcal{E}_\lambda, \widehat{\mathcal{E}}_\lambda) \stackrel{a.s.}{\to} 0.$$

*Hence, for every $\xi > 0$, there exists $\lambda > 0$ such that $d_H(\mathcal{A}_\sigma, \widehat{\mathcal{E}}_\lambda) = O_P(\sigma + \xi).$*

REMARK 5. The theorem shows that the level sets approximate $\mathcal{A}_\sigma$ which is itself and approximation to $\mathcal{A}$. We know that $\mathcal{A}_\sigma$ is close to $\mathcal{A}$ but a more precise statement would require specific assumptions on the shape of the filaments.

4.3. *Examples.* In this section we briefly consider an example based on galaxy data and a synthetic example.

EXAMPLE 1. For the simulated example, the vertices of a pentagon have been randomly selected over $[0, 1]^2$. On each side of the pentagon, points were drawn according to a Beta$(\frac{1}{2}, \frac{1}{2})$ distribution. The total of $n = 500$ points is divided among the sides proportionally to the side lengths. The resulting data set, shown in Figure 2A, was obtained adding a bivariate normal perturbation (with $\sigma = 0.03$) to the points.

Figure 2B and C shows that all the steepest ascent paths move toward the perimeter of the pentagon and then along it until they reach a mode. The level set at the 90th percentile of the estimated path density seems to be a good approximation to the pentagon perimeter (see Figure 2D).

To see the effect of noisy background, 500 uniform points in $[0,1]^2$ were added to the data set. The whole procedure was implemented again on the augmented data set. Figure 8 shows that the presence of background introduces some noise, but the procedure still catches the filaments.

EXAMPLE 2. The filament estimation procedure has been tested also on cosmological data. The data set analyzed is part of the mini "SDSS ugriz" catalogue, produced by Croton et al. (2006) and publicly available online at www.mpa-garching.mpg.de/galform/agnpaper/. Along with other information, the catalogue gives the position of all galaxies in the cubic box $[0, 62.5 \text{ Mpc/h}]^3$. To deal with a two-dimensional data set, we selected the galaxies whose $z$-coordinate is in the interval $[20, 25]$ Mpc/h and considered only the first two coordinates. The resulting data set, shown in Figure 1A, contains the $(x, y)$-coordinates of $n = 2435$ galaxies. Plots B, C and D in Figure 1 show a reasonably satisfactory behavior of the filament detection procedure.

**5. Discussion.** Our results are a first step toward developing methods for finding filaments without imposing strong assumptions about those filaments. However, much work needs to be done to turn the theory into practical methodology.

First, we need good, data-driven methods for choosing the bandwidths $h_n$ and $\nu_n$. We conjecture that cross-validation methods could be quite effective.

Second, extracting the filaments from the path density estimator $\widehat{p}_n$ is obviously very important. To find high local concentrations of paths, one can trim the observed paths or choose high level sets, as we did in the examples. These approaches, and possibly others, deserve careful examination.

Third, we conjecture that these methods will generalize to three dimensions. The main challenge in doing so lies in the characterization of the path density, which requires handling an additional class of singular points and generalizing part of the proof from using closed curves to using closed surfaces. Both steps are straightforward but tedious. Generalization to higher dimensions remains an open but interesting problem; related problems in computational topology have proven difficult in high dimensions in some cases.

Finally, we think that the mathematical methods we have developed could be useful in other problems as well.

**6. Proofs for Section 2.** We begin with an explicit characterization of the path density $p$ that includes the statement of Theorem 1. Throughout, we let $\mathcal{U} = \mathbb{R}^2$ and let $U_0 \subset \mathcal{U}$ be a large compact set that contains all the critical points of $g_X$ in its interior.

THEOREM 6. *Under Assumptions 1–3, the path density $p$ is an upper-semicontinuous function with the following properties:*

1. *If $x$ is a local minimum of $g_X$, $p(x) = 0$.*
2. *If $x$ is a local maximum of $g_X$, $p(x) = \infty$.*
3. *If $x$ is a saddlepoint of $g_X$, then $0 < p(x) < \infty$ and there exist four sequences of regular points $(x_k^j)$, for $j = 1, \ldots, 4$, all converging to $x$ such that*

$$p(x) = \lim_{n \to \infty} \sum_{j=1}^{4} p(x_n^j). \tag{18}$$

4. *If $x$ is a regular point of $g_X$, then the following hold:*
   (a)
   $$p(x) = \lim_{r \to 0} \frac{\pi(\partial B(x,r))}{r}. \tag{19}$$

   (b) *For all sufficiently small $r > 0$, there exists a curve $\gamma$ of and $\alpha_r < 0 < \beta_r$ such that* (i) $\text{length}(\gamma[\alpha_r, 0]) = r + O(r^2)$, (ii) $\text{length}(\gamma[0, \beta_r]) = r + O(r^2)$ *and* (iii)
   $$p(x) = \lim_{r \to 0} \frac{\pi(\gamma([\alpha_r, \beta_r]))}{r}. \tag{20}$$

   (c) *There exists a smooth function $h : \mathcal{U} \times \mathbb{R} \to [0, \infty)$ such that for any regular point $x$ of $g_X$*
   $$p(\psi^{-t}x) = \int_0^\infty g_X(\psi^{-t-s}x) h(x, t+s) \, ds. \tag{21}$$

5. *$p$ is continuous on $(\mathcal{M} \cup \mathcal{H})^c$ and bounded on $\mathcal{M}^c$.*

To prove Theorem 6 we need two lemmas whose proofs are reported at the end of this section. Let $D^2 g_X(x)$ denote the Hessian matrix of $g_X$ at $x$.

LEMMA 1. *For $x \in U_0$, let $r_0 > 0$ be small enough so that $B(x, r_0)$ contains no critical points of $g_X$, other than possibly $x$ itself. Then, for $r \leq r_0$, there exists a constant $C > 0$ such that for $y \in B(x, r)$*
$$|g_X(y) - g_X(x) - \nabla g_X(x) \cdot (y - x) - \tfrac{1}{2}(y-x)^T D^2 g_X(x)(y-x)| \leq C r^3,$$
$$\|\nabla g_X(y) - \nabla g_X(x) - D^2 g_X(x)(y-x)\| \leq C r^2,$$
*where $C$ is independent of $r$ and $y$ but possibly dependent on $r_0$.*

LEMMA 2. *Let $\gamma:[-1,1] \to U_0$ be a nonintersecting curve such that $\gamma([-t,t])$ has length $\ell_t$ and such that $\gamma(t)$ is a regular point of $g_X$ for every $-1 \leq t \leq 1$. Then, there exists a $0 < t \leq 1$ such that, for all $0 < s < t$,*

$$(22) \qquad \mathcal{L}(\mathcal{V}(\gamma([-s,s]))) \leq C\ell_s + O(\ell_s^2),$$

*where $C > 0$ is independent of $s$ but possibly dependent on $x$ and $t$.*

PROOF OF THEOREM 6. 1. The idea of the proof is to show that in small enough balls around $x$, the integral curves closely approximate those from a quadratic function around $x$. If $g_X$ were quadratic at $x$, then the flows $\psi^{-t}$ carry every point in $B(x,r)$ radially toward $x$. Hence, $\pi(B(x,r)) \leq \mu_X(B(x,r)) = O(r^2)$. The details that follow require some careful arguments because we make minimal assumptions about $g_X$ other than smoothness.

We begin by constructing two small open balls $W_0 \subset V$ around $x$ with four properties: (i) $V$ contains no critical points of $g_X$ other than $x$, (ii) the Laplacian $\nabla^2 g_X > 0$ on $V$, (iii) $g_X$ is locally quadratic on $V$ in the sense that $g_X(y) = g_X(x) + h_1^2(y) + h_2^2(y)$, for a smooth function $h$, and (iv) the reverse evolution of $W_0$ remains within $V$, $\overline{\mathcal{V}(W_0)} \subset V$.

Note first that $\psi^t x = x$ for all $t \in \mathbb{R}$. Because the Laplacian $\nabla^2 g_X$ is continuous and positive at local minimum $x$, there is an open ball $V_1$ around $x$ on which $\nabla^2 g_X > 0$. By Morse's theorem [Milnor (1963)], there is a neighborhood $V_2$ of $x$ and a $C^\zeta$ diffeomorphism $h$ from $V_2$ to an open ball $A$ around 0 in $\mathbb{R}^2$ such that $h(x) = 0$ and for $y \in V_2$,

$$(23) \qquad g_X(y) = g_X(x) + h_1^2(y) + h_2^2(y).$$

Note $V_2$ satisfies (i). Let $V$ be an open ball around $x$ contained in $V_1 \cap V_2$. Then $V$ satisifies properties (i), (ii) and (iii) above.

To find $W_0$ satisfying (iv) we need to get close enough so that the quadratic behavior dominates any wandering of the integral curves outside of $V$. Define $q: A \to \mathbb{R}$ by $q(z) = \|h^{-1}(z) - x\|^2$. This is a $C^\zeta$ function with $q(0) = 0$ and hence Lipschitz. Thus, there exists an open ball $A_1 \subset A$ such that $q(z) \leq c\|z\|$ on $A_1$, for some constant $c$. Hence, there is an open ball $W_0 \subset h^{-1}(A_1)$ around $x$ such that $\overline{\mathcal{V}(W_0)} \subset V$.

Next we show that for any open ball $W \subset W_0$ around $x$, its reverse evolution must stay within its closure, $\mathcal{V}(W) \subset \overline{W}$. This is equivalent to showing that integral curves from within $W$ cannot cross the boundary $\partial W$. To do this, we will use the following result. Let $C \subset V$ be any region bounded by a piecewise smooth closed curve. By Stokes's theorem and by property (ii) of $V$ above,

$$(24) \qquad \int_{\partial C} (-\nabla g_X) \cdot \hat{\mathbf{n}} = \int_C (-\nabla^2 g_X) < 0,$$

where $\hat{\mathbf{n}}$ is an outward pointing, unit vector normal to $\partial C$. (The application of Stokes's theorem comes by defining the differential form $\omega = -v_2\, dx + v_1\, dy$, where $v = -\nabla g_y$. Then, $d\omega = -\nabla^2 g_X\, dx \wedge dy$ and $\int_C \omega = \int v \cdot \hat{\mathbf{n}}$.)

Apply (24) with $C = W$. If the integrand on the left, $(-\nabla g_X) \cdot \hat{\mathbf{n}}$, is ever positive on $\partial W$, it must, by continuity, be positive on some arc $(a_1, a_2) \subset \partial W$. Suppose then that the integrand on the left is nonnegative on such an arc $(a_1, a_2)$. We will construct a closed curve bounding a region $C \subset V$ to which we will apply (24). Construct the closed curve by jointing the following three curves, in order: a curve from $x$ to $a_1$ obtained by $t \in [0,1] \mapsto h^{-1}(th(a_1))$, the arc $(a_1, a_2)$ and the curve from $a_2$ to $x$ obtained by $t \in [0,1] \mapsto h^{-1}((1-t)h(a_2))$. Note that the first and last of these curves have two important properties: (i) their image is contained in $V$, and (ii) they follow integral curves of $\psi$ in one direction or the other. The latter is implied directly by (23). As a result, the integrand on the left-hand side of (24) is zero over the first and last segments and nonnegative over the middle one, which is a contradiction equation (24).

As a result, for sufficiently small $r > 0$, an open ball $B(x,r)$ satisfies $\mathcal{V}(B(x,r)) \subset \overline{B(x,r)}$, and thus

$$(25) \qquad p(x) = \lim_{r \to 0} \frac{\pi(B(x,r))}{r} \leq \lim_{r \to 0} \frac{O(r^2)}{r} = 0.$$

2. We can apply the same argument as in item 1 except with Morse's theorem giving the representation

$$(26) \qquad g_X(y) = g_X(x) - h_1^2(y) - h_2^2(y)$$

and thus reversing the sign of the Laplacian. This shows that for sufficiently small $r_0 > 0$, an open ball $B(x,r_0)$ satisfies $\mathcal{V}(B(x,r_0)) \supset B(x,r_0)$. Moreover, this implies that for any $r < r_0$, $\mathcal{V}(B(x,r)) \supset B(x,r_0)$. Hence,

$$(27) \qquad p(x) = \lim_{r \to 0} \frac{\pi(B(x,r))}{r} \geq \lim_{r \to 0} \frac{\mu_X(B(x,r_0))}{r} = \infty.$$

3. Note again that $\psi^t x = x$ for all $t \in \mathbb{R}$ and again that Morse's theorem implies the existence of a neighborhood $U$ of $x$ and a diffeomorphism $h$ from $U$ to an open ball $A$ around 0 in $\mathbb{R}^2$ such that $h(x) = 0$ and for $y \in U$,

$$(28) \qquad g_X(y) = g_X(x) - h_1^2(y) + h_2^2(y).$$

Define $\xi$ as the global flow generated by the vector field $V^\perp = R\nabla g_X$, where $R$ is a fixed 90-degree rotation matrix.

Choose $r > 0$ such that $B(x,r) \subset U$ and $B(x,r)$ contains no critical points other than $x$. For every $0 < t \leq r$, let $S_t = h(\partial B(x,t))$. This is a smooth closed curve around 0 in $A$. In the $h$ coordinate system in $A$, the integral curves of the vector field are hyperbolas of the form $h_1 h_2 = c$ for some constant $c$, with the reversed (negative gradient) flow traveling from large $|h_2|$ and small $|h_1|$ to small $|h_2|$ and large $|h_1|$. (Here we are treating $h_1$ and $h_2$ as the coordinates under the change of variables induced by the diffeomorphism $h$.) Specifically, in $A$, the flow through a point $(h_1, h_2)$ takes the parametric form $s \mapsto (h_1 e^{2s}, h_2 e^{-2s})$. It follows that there is at least one point in each quadrant where the product $|h_1 h_2|$ attains its maximum on $S_t$, and that the set of such points in each quadrant is a closed set (by continuity). It follows that if there is more than one such point, they must all be on the same contour and that therefore, because the set is compact, there must be a unique point that succeeds the others under the partial order $\succeq$. Let $z_i(t) \in S_t$, for $i = 1, \ldots, 4$, denote these points, and let $y_i(t) = h(z_i) \in \partial B(x, t)$. In addition, let $z_+(t)$ and $z_-(t)$ denote the points at which $S_t$ intersects the positive and negative $h_1$ axis, respectively, and let $y_+(t) = h(z_+(t))$ and $y_-(t) = h(z_-(t))$. Here, the subscripts $1, \ldots, 4$ indicate the associated quadrant and the subscripts $+$ and $-$ indicate association with the positive and negative $h_1$ axis, respectively. Take $z(0) = 0$ and $y(0) = x$ in each case.

The six curves in the original space $y_1, y_2, y_3, y_4, y_-, y_+$ are smooth and parameterized by arclength. Moreover, since the integral curves are tangent at each of the points $y_i(r)$, for $i = 1, \ldots, 4$, each curve $y_i$ traces an integral curve of the flow $\xi$.

From these, we construct four closed curves, $\gamma_1, \gamma_2, \gamma_3$ and $\gamma_4$, one per quadrant, as follows. The first travels along $y_1$ from $x$ to $y_1(r)$, then along the arc of $\partial B(x, r)$ to $y_+(r)$, then along $y_+$ from $y_+(r)$ back to $x$. The rest are analogous, traveling $x \to y_+(r) \to y_2(r) \to x$, $x \to y_4(r) \to y_-(r) \to x$ and $x \to y_-(r) \to y_3(r) \to x$. We choose the time index along the circular arcs so that the $\gamma_i$'s are parameterized by arclength. Note that the range of $\gamma_i$ is restricted to quadrant $i$ by construction.

Let $C$ denote the region (including the boundary) enclosed by $\gamma$ curves. This consists of two "lobes," one joining the region enclosed by the two $\gamma$ curves that intersect along the range of $y_+$ and the other joining the regions enclosed by the two $\gamma$ curves that intersect along the range of $y_-$. The two lobes intersect only at $x$. Partition $B(x, r) = Z + L + R$, where $Z = B(x, r) \cap h^{-1}(\{(u, 0) \in A : u \neq 0\})$ is the image of the $h_1$ axis; $L = B(x, r) \cap (C - Z)$ is the union of the two lobes minus the axis; and $R = B(x, r) - L - Z$ is everything else. The construction of the $y_i$'s shows that these sets satisfy the following: (i) every point in $R$ preceeds a point [which is neither $x$ nor the $y_i(r)$'s] on one of the $y_i$ curves; (ii) $Z$ is the image of a curve (and thus a set of measure 0), so every point in $Z$ precedes (under $\preceq$) one of the two points in $\{y_+(r), y_-(r)\}$; (iii) every point in $L$ is succeded (under $\succeq$) by a point (which is neither $x$ nor the $y_i$'s) on one $y_i$ piece of the corresponding curve. Note that the integral curves from each $y_i$ piece do not cross the "axis" $h(\{(u, 0) \in A : u \neq 0\})$. (These facts are easiest to see in the transformed space $A$, even though the curves need not be radial or arc-like in that space.

The $z_i$'s are the latest points on the highest integral curves. $Z$ is the image of the $h_1$ axis which creates two integral curves.)

It follows then that

(29) $\quad \mathcal{V}(B(x,r)) = \{x\} \cup \mathcal{V}(\gamma_1((0,r])) \cup \cdots \cup \mathcal{V}(\gamma_4((0,r])) \cup R \cup \mathcal{V}(Z),$

where this union is disjoint. Moreover, $\mathcal{L}(R) \leq \pi r^2$ and $\mathcal{L}(\mathcal{V}(Z)) = 0$, so neither component contributes to $p(x)$. Because the $y_i$'s follow the integral curves of $\xi$, for any $0 < t < r$, $p(y_i(t)) = \lim_{\delta \to 0} \pi(y_i((t-\delta, t+\delta)))/\delta$ by Theorem 6(4b)—proved below, independently of these results—and for disjoint intervals along $y_i$, the corresponding reverse evolutions are disjoint. For integer $n > 0$, let $\Delta = r/n$, then,

$$\pi(y_i((0,r))) = \sum_{j=1}^{n} \pi(y_i((j-1)\Delta, j\Delta))$$

(30)
$$= \sum_{j=1}^{n} \frac{\pi(y_i((j-1)\Delta, j\Delta))}{\Delta} \Delta$$

$$\to \int_0^r p(y_i(s)) \, ds$$

as $n \to \infty$. The convergence is justified by the dominated convergence theorem as each of the functions $f_n(t) = \frac{\pi(y_i(t-\Delta/2, t+\Delta/2))}{\Delta}$ converges pointwise to the bounded function $p(y(t))$. It follows then that

$$\lim_{r \to 0} \frac{\pi(y_i((0,r)))}{r} = \lim_{r \to 0} \frac{1}{r} \int_0^r p(y_i(s)) \, ds = \lim_{r \to 0} p(y_i(r)).$$

And therefore,

(31) $$p(x) = \sum_{i=1}^{4} \lim_{r \to 0} p(y_i(r)),$$

which proves the claim.

For the upper semicontinuity proof below, note that any sequence $x_n \to x$ is eventually all regular points, and each $x_n$ can lie in only one "quadrant" (relative to $A$), so $\limsup_n p(x_n)$ is no greater than the maximum of the lim sup over any subsequence lying in one quadrant. Hence, $p(x) \geq \limsup_{n \to \infty} p(x_n)$, so $p$ is upper-semicontinuous at $x$.

(4a) For small enough $r$, $B(x,r)$ contains no critical points of $g_X$. Partition $\partial B(x,r)$ into two sets $C_{\text{in}}$ and $C_{\text{out}}$, where $C_{\text{out}}$ contains all points on $\partial B(x,r)$, whose paths do not cross the interior of $B(x,r)$. Every point in $C_{\text{in}} = \partial B(x,r) - C_{\text{out}}$ either precedes or succeeds a point in $B(x,r)$. So, $\mathcal{V}(C_{\text{in}}) \subset \mathcal{V}(B(x,r)) \cup \partial B(x,r)$. Because the flow is smooth with a bounded derivative on $U_0$ and because $\partial B(x,r)$ is compact, $C_{\text{out}}$ equals an at most countable union of points and arcs on which the flow is tangent to the circle. Each such piece produces a single curve under $\mathcal{V}$, so it follows that $\mathcal{L}(\mathcal{V}(C_{\text{out}})) = 0$, which in turn implies that $\pi(C_{\text{out}}) = 0$. Because the symmetric difference of $\mathcal{V}(B(x,r))$ and $\mathcal{V}(\partial B(x,r))$ is contained in $\overline{B(x,r)} \cup \mathcal{V}(C_{\text{out}})$, we have that $\pi(B(x,r)) = \pi(\partial B(x,r)) + O(r^2)$ as $r \to 0$.

(4b) For convenience, let $u = \nabla g_X : \mathbb{R}^2 \to \mathbb{R}^2$ denote the gradient function. Let $u^0 = u(x)$ and let $H$ be the Hessian of $g_X$, with eigenvalues $\lambda_1 > \lambda_2$. Pick $\delta > 0$ small. Choose $r_0 > 0$ such that:

(i) $B(x, r_0)$ contains no critical points of $g_X$;
(ii) $r_0 \frac{\lambda_1 - \lambda_2}{2\|u^0\|} \leq 1/4$;
(iii) $r_0 \frac{\lambda_1 + \lambda_2}{2\|u^0\|} \leq 3/4$;
(iv) $\|u(y) - u(x) - H(y-x)\| \leq \delta(\|y-x\|/r_0)^2$ via Lemma 1.

Then, let $r < r_0$.

The plan of the proof is to construct the curve $\gamma$ explicitly in the case where $u$ is linear and show that the perturbation caused by the addition of higher order terms causes only a small change to the curve, consistent with the statement of the theorem. The range of the curve will generate the same set under $\mathcal{V}$ as the open ball around $x$, up to an $O(r^2)$ term. Consider points along the circles around $x$ where the gradient $u$ is tangent to the circle. Connecting these points will cut all the integral curves within the ball. Note that because the vector field $\nabla g_y$ is curl free, it follows from Stokes' theorem that there must exist at least two tangent points. That is, because the line integral around $S$ is zero, there must be a sign change of the tangent vectors, but this requires at least two zeros on the circle. [The application of Stokes' theorem comes by defining the differential form $\omega = v_1 \, dx + v_2 \, dy$, where $v = \nabla g_y$. This form is closed, that is, $d\omega = (D_1 v_2 - D_2 v_1) \, dx \wedge dy = 0$, and so $0 = \int_C \omega = \int v(C(t)) \cdot \frac{dC}{dt}(t)$.]

Let $B_0 = \overline{B(x, r_0)}$, $B = B(x, r)$ and $S = \partial B$. We begin by showing that when $u$ is linear, it is tangent to $S$ at exactly two isolated points for every $r < r_0$. We show further that the component of $u$ normal to $S$ is small only in a small neighborhood of these tangent points, which will be used to show that such tangent points lie on two small arcs for general $u$.

So, for the moment, assume that $u(y) = u(x) + H(y-x)$ on $B_0$. Recall that $u \neq 0$ on $B_0$. Because $H$ is symmetric, there are two orthogonal, unit vectors $v_1$ and $v_2$ and two real numbers $\lambda_1 \geq \lambda_2$ (possibly zero) such that $H(\beta_1 v_1 + \beta_2 v_2) = \beta_1 \lambda_1 v_1 + \beta_2 \lambda_2 v_2$. For every $y \in B$, $y - x$ can be written as $\beta_1 v_1 + \beta_2 v_2$ for $\beta_1^2 + \beta_2^2 \leq r^2$ and vice versa. Hence, we can write, with some abuse of notation,

$$(32) \qquad u(\beta_1, \beta_2) = u(0,0) + \lambda_1 \beta_1 v_1 + \lambda_2 \beta_2 v_2,$$

where $u(\beta_1, \beta_2) = u(x + \beta_1 v_1 + \beta_2 v_2)$. Let $\alpha_1, \alpha_2$ be such that $u(0,0) \equiv u(x) = \alpha_1 v_1 + \alpha_2 v_2$. We have by assumption of regularity at $x$ that $\|\alpha\| = \sqrt{\alpha_1^2 + \alpha_2^2} > 0$. Then, any $y \in S$ corresponds to $\beta = (r\cos\theta, r\sin\theta)$ for some $\theta \in [0, 2\pi)$, and

$$(33) \qquad u(\beta_1, \beta_2) = (\alpha_1 + \lambda_1 r \cos\theta) v_1 + (\alpha_2 + \lambda_2 r \sin\theta) v_2.$$

Because $u$ is nonzero on $B_0$, we have that for $y \in S$, $u(y)$ is tangent to $S$ if $u(\beta_1, \beta_2) \cdot (\beta_1, \beta_2) = 0$, or more explicitly,

$$(34) \begin{aligned} 0 &= [(\alpha_1 + \lambda_1 r \cos\theta) v_1 + (\alpha_2 + \lambda_2 r \sin\theta) v_2]^T (r\cos(\theta) v_1 + r\sin(\theta) v_2) \\ &= r(\alpha_1 \cos\theta + \lambda_1 r \cos^2\theta) + r(\alpha_2 \sin\theta + \lambda_2 r \sin^2\theta), \end{aligned}$$

which implies, canceling $r > 0$, that

$$(35) \qquad \alpha_1 \cos\theta + \alpha_2 \sin\theta + r\left(\frac{\lambda_1 - \lambda_2}{2}\cos 2\theta + \frac{\lambda_1 + \lambda_2}{2}\right) = 0.$$

Because $\|u^0\| = \|u(x)\| = \|\alpha\| > 0$, there exists $\eta \in [0, 2\pi)$ such that $\cos\eta = \alpha_1/\|\alpha\|$ and $\sin\eta = \alpha_2/\|\alpha\|$. Then, solutions of the above equation correspond to crossings of the purely imaginary axis of the following complex curve $c:[0, 2\pi) \to \mathbb{C}$:

$$(36) \qquad c(\theta) = e^{i(\theta-\eta)} + r\frac{\lambda_1 - \lambda_2}{2\|\alpha\|}e^{i2\theta} + r\frac{\lambda_1 + \lambda_2}{2\|\alpha\|}.$$

This is a generalized epicycloid with phase $\eta$ and offset $r\frac{\lambda_1+\lambda_2}{2\|\alpha\|}$ along the real axis. Write $c_1 = \text{Re}(c)$, $c_2 = \text{Im}(c)$, $w = r\frac{\lambda_1-\lambda_2}{2\|u^0\|} \geq 0$ and $v = r\frac{\lambda_1+\lambda_2}{2\|u^0\|}$. Following Brannen (2001) and treating $(c_1 - v, c_2)$ as a plane curve, we see that the curve has no cusps if the curvature never changes sign and no loops if the vector cross-product of the curves position and velocity never changes sign. The sign of the curvature equals the sign of

$$(37) \qquad c_1' c_2'' - c_1'' c_2' = 8w^2 + 1 + 6w\cos(\theta + \eta)$$

and the sign of the vector cross-product equals the sign of

$$(38) \qquad (c_1 - v)c_2' - c_1' c_2 = 2w^2 + 1 + 2w\cos(\theta + \eta).$$

Taking $0 \leq w \leq 1$ and using the fact that the cosine terms are between $-1$ and $1$, we see that the former keeps the same sign if $w \leq 1/4$ or $w \geq 1/2$, and the latter keeps the same sign if $(2w-1)(w-1) \geq 0$, which requires $w \leq 1/2$. The curve is thus locally convex whenever $w \leq 1/4$. Moreover, the $(c_1 - v)^2 + c_2^2 \geq (1-w)^2$, so the curve always lies outside a circle of radius $1-w$ around $v$. Hence, if $0 \leq w \leq 1/4$ and $|v| < 1-w$, the curve will intersect the imaginary axis exactly twice. By assumption, then, the original vector field thus has exactly two isolated points of tangency with each circle $S$.

By the implicit function theorem, the collection of these tangent points, together with $x$, form the image of a smooth curve, $\gamma$, with $\gamma(0) = x$ and $\gamma'(0)$ nonzero. Every point along this curve is perpendicular to the gradient flow, so this curve is part of the integral curve through $x$ for the flow $\xi$, so we can take $\gamma(t)$ of the form $\xi^t x$. The angular separation between the two tangent points forming the branches of $\gamma$ can be determined by studying the curve $c$ above. Direct calculation shows that the arc length varies with $\theta$ along the curve as $1 + O(r^2)$ relative to a circle. It follows that the angle between the tangent points is $\pi + O(r)$, and thus the length of $\gamma$ between the center of the circle and each point on the circle at radius $r$ has length $r + O(r^2)$. This proves the claim for the linear case.

In the general case, an error term of order $r^2$ is added to linear term by the approximation Lemma 1. This perturbs the tangent points and can add or move them, though as shown above there are always at least two. (We get them by following $\xi$ forward and backward until we leave the circle, which we do eventually because there are no critical points in the neighborhood.) We show that these are confined to two small arcs on the circle of length $O(r^4)$. To do this, note that the size of the outward normal component of $u(y)$ is $|u(y) \cdot (y-x)/\|y-x\|| = |u(y) \cdot (y-x)|/r$. The numerator is the left-hand side of (35). The normal component might be zero for any $\theta \in [0, 2\pi)$ for which $|\operatorname{Re}(c(\theta)/r)| \leq Cr^2$, or equivalently $|\operatorname{Re}(c(\theta))| \leq Cr^3$. For a circle, this occurs over two arcs of angular size $O(r^3)$ and thus of length $O(r^4)$ in $U_0$. For the general epicycloid, the arc length over any segment can be at most $1 + 2w$ times that of the circle, which again gives an arc length of $O(r^4)$.

The $\gamma$ curve constructed above connects two tangent points, the remaining integral curves can cross the boundary in a curve of length $O(r^4)$ and hence by Lemma 2,

$$(39) \qquad \mathcal{L}(\mathcal{V}(B(x,r)) - \mathcal{V}(\gamma([\alpha_r, \beta_r]))) = O(r^4),$$

where $\alpha_r$ and $\beta_r$ are the time indices at which each branch of $\gamma$ strikes the circle of radius $r$. The result follows.

It is worth noting that we can construct a curve $\widetilde{\gamma}$ for which $\mathcal{V}(B(x,r))$ equals $\mathcal{V}$ applied to the image of the curve, except for part of $B(x,r)$ itself. To do this, we include all the tangent points on the two arcs of length $O(r^4)$ and then connect opposite ends with a curve of minimal length through $x$. The resulting curve has length $2r + O(r^2)$ and hits every equivalence class in $B(x,r)$.

(4c) By Theorem 6(4b) above, it is sufficient to work with $\mathcal{V}(\gamma([a_r, b_r]))$ for sufficiently small $r > 0$. This is generated by a curve of the form $\gamma(s) = \xi^s x$ over the interval $[a_x(r), b_x(r)]$. Taking $a_x(0) = b_x(0) = 0$, note that $a_x$ and $b_x$ are differentiable functions of $r$ and continuous functions of $x$. Hence, we have the set $\mathcal{A} = \mathcal{V}(\gamma([a_x(r), b_x(r)]))$. Construct a smooth bijection $[0, \infty) \times [a(r), b(r)]$ to $\mathcal{A}$ by $z(t, s) = \psi^{-t}\xi^s x$. It follows by the change of variables and Fubini theorems that

$$
\begin{aligned}
\frac{\pi(\mathcal{V}(\gamma_r([a_r, b_r])))}{r} &= \frac{1}{r} \int_A g_X(v)\, dv + O(r^2) \\
&= \frac{1}{r} \int_{a(r)}^{b(r)} \int_0^\infty g_X(\psi^{-t}\xi^s x) J(x, t, s)\, dt\, ds + O(r^2),
\end{aligned}
\tag{40}
$$

where $J(x, t, s)$ is the Jacobian determinant of $z$, which is a smooth function in $x$. Taking limits as $r \to 0$ yields

$$
p(x) = (b'_x(0) - a'_x(0)) \int_0^\infty g_X(\psi^{-t}\xi^0 x) J(x, t, 0)\, dt. \tag{41}
$$

Taking $h(x, t) = (b'_x(0) - a'_x(0)) J(x, t, 0)$ gives the theorem for $p(x)$. Now apply the above to $\psi^{-t} x$ for any $t \in \mathbb{R}$ to get the formal statement.

Boundedness of $p$ on $U - \mathcal{M}$ follows from showing uniform boundedness at regular points of $g_X$, by items 1 and 3 above. At a regular point $x$ of $g_X$ and for sufficiently small $r > 0$, there is by Theorem 6(4b) above and Lemma 2 a constant $K > 0$ dependent only on $g_X$ such that $\pi(\gamma([a_x(r), b_x(r)])) \leq Kr$. It follows that $p(x) \leq K$ as was to be proved.

Finally, to prove that $p$ is upper semi-continuous it is sufficient to show that $p$ is continuous at regular points and local minima. The rest follows from items 1 and 2, and the proof of item 3 above.

Suppose $x$ is a local minimum of $g_X$ and pick $a_0 > 0$ small. Then by the proof of item 1 above, there exists an $a_1 \leq a_0$ such that $a \leq a_1$ implies that $\mathcal{V}(B(x, a)) \subset \overline{B(x, a)}$. It follows from Lemma 2 that if $y \in B(x, a)$, then for small $r < a/2$, $\mathcal{L}(\mathcal{V}(B(y, r))) \leq Cra$, the width of the ball times the scale of the set $B(x, a)$, where $C > 0$ is a constant independent of $a$ but dependent possibly on $a_0$. If we choose $a < a_0/C$, then $\pi(B(y, r))/r \leq a_0$ for $r < a/2$, showing that $p(y) \leq a_0$. It follows that $p(y) \to 0$ as $y \to x$, which proves continuity of $p$ at local minima.

To prove continuity at regular points, we apply Theorem 6(4c) above. Because $g_X$ and the $h$ in Theorem 6(4c) are continuous in $x$, they are bounded on $U_0$. The bounded convergence theorem thus allows an interchange of limit and integral. Taking $t = 0$ and letting $x_n \to x$ all be regular points of $g_X$, we thus have that

$$
p(x) = \lim_{n \to \infty} \int_0^\infty g_X(\psi^{-s} x_n) h(x_n, s)\, ds = \lim_{n \to \infty} p(x_n). \tag{42}
$$

This proves that $p$ is upper semi-continuous. This completes the proof. $\square$

PROOF OF THEOREM 2.   1. By Theorem 6(4b) and (c), we can write

$$
\pi(B(x, r)) = \int_{\alpha_r}^{\beta_r} \int_0^\infty H(\xi^s x, t)\, dt\, ds + O(r^2), \tag{43}
$$

where by the proof of Theorem 6(4c), $H$ is a smooth function and $\beta_r - \alpha_r = r + O(r^2)$. If $D_1 H(y, t)$ denotes the derivative of $H$ with respect to its first argument, then Taylor's theorem gives us that $H(y, t) = H(x, t) + D_1 H(u(t, y), t) \cdot (y - x)$ for some points $u(t, y)$ on the line between $y$ and $x$. Because $H(y, t)$ and $H(x, t)$ are integrable with respect to $t$, and $y$ is arbitrary in $\mathcal{U}$, it follows that $\|D_1 H(u(t, y), t)\|$ is integrable as well. The integrals are bounded over $U_0$ by compactness. Hence, for $s \in [\alpha_r, \beta_r]$,

$$(44) \qquad \left| \int_0^\infty H(\xi^s x, t)\, dt - \int_0^\infty H(x, t)\, dt \right| \le Cr$$

for a constant $C > 0$ independent of $x$. The second term in the above absolute value is just $p(x)$, so

$$(45) \qquad \pi(B(x, r)) = \int_{\alpha_r}^{\beta_r} (p(x) + O(r))\, ds + O(r^2) = rp(x) + O(r^2),$$

which proves the claim.

2. Let $\gamma$ be the curve in Theorem 6(4b) centered at $x$ and parameterized by arc-length over the interval $[a, b]$, where $a < 0 < b$, $|a| = r + O(r^2)$, and $|b| = r + O(r^2)$. For positive integer $m$, pick index points $t_{mk} = a + (b - a)k/m$, for $0 \le k \le m$ and cover $\gamma([a, b])$ by open balls of radius $(b - a)/2m$ centered at each $\gamma(t_{mk})$. It follows that

$$(46) \qquad \pi(\gamma([a, b])) = \sum_{k=0}^m \pi\left( B\left( \gamma(t_{mk}), \frac{1}{2m} \right) \right) + O(m^{-2})$$

$$(47) \qquad \approx \sum_{k=0}^m p(t_{mk}) \frac{1}{2m} + O(m^{-2})$$

$$(48) \qquad \to \frac{1}{2} \int_a^b p(\gamma(t))\, dt$$

$$(49) \qquad = \frac{1}{2} \int_{-r}^r p(\gamma(t))\, dt + O(r^2),$$

where the limit is over $m \to \infty$. The second line above follows by part 1. Because $\pi(B(x, r)) = \pi(\gamma([a, b])) + O(r^2)$, the result follows. □

PROOF OF LEMMA 1. This result follows directly from Taylor's theorem and the compactness of $U_0$. For example, the remainder term in the gradient approximation (second equation in the theorem) can be written $r^2 (u^T A_1(y, x) u,\ u^T A_2(y, x) u)^T$, where $u = (y - x)/\|y - x\|$ and where the $A_i(y, x)_{jk} = \int_0^1 (1 - t) \partial^3 g_X(x + t(y - x))/\partial x_i\, \partial x_j\, \partial x_k\, dt$. Because the maximum eigenvalue of these matrices is a continuous function of the matrix and thus continuous in $y$ [Naulin and Pabst (1994)], each component of this vector is bounded. □

PROOF OF LEMMA 2. To begin, suppose $\gamma$ is a short segment of length $r$ of a curve of the form $\gamma(s) = \xi^s x$, parameterized on the interval $[a,b]$ with $a < 0 < b$. Assume that every point along the curve is a regular point of $g_X$. For sufficiently small $r$, the indices $a < 0 < b$ satisfy $|b - a| \leq \kappa r$ because the derivative of $\nabla g_X$ is bounded above and below (in norm) in a small neighborhood of $x$.

Every point $\gamma(s)$ can be classified according to whether $\lim_{t \to \infty} \psi^{-t} \gamma(s)$ is (i) a point on the boundary of $U_0$; (ii) a local minimum of $g_X$; or (iii) a saddle point of $g_X$. If there are no points of class (iii) on the segment, then we can take $r$ small enough so that all points are of the same class. In general, under Assumption 2, there is at most a finite set of $s \in [a,b]$ of class (iii), each of whose $\mathcal{V}(\gamma(s))$ has Lebesgue measure 0. So it suffices to assume that there are no class (iii) points on $[a,b]$ because without loss of generality, we could individually consider the finite open intervals in $[a,b]$ on which this is the case.

Let $V$ denote the set $\mathcal{V}(\gamma([a,b]))$ and let $I = \bigcup_{a \leq s \leq b} I_{\gamma(s)}$, where $I_x$ is the intervals of time indices for which the flow from $x$ stays in $U_0$. Define

$$V_1 = \{\psi^{-t} \xi^s x : s \in [a,b], t \in I\} \subset \mathcal{U}$$

and note that $V \subset V_1$. Define a mapping $h : [0, \ell] \times I \to V_1$ by $h(s,t) = \psi^{-t} \xi^s x$. Note that if $s_1 \neq s_2$, then the $\xi^{s_i} x$ are in different equivalence classes and so $h(s_1, t) \neq h(s_2, t)$ for any $t$. Similarly, if $t_1 \neq t_2$, then $\psi^{-t_i} \xi^s x$ lie at different points along the integral curve. Hence, $h$ is one to one. If $y \in V_1$, then by construction, $y$ is equivalent to a unique $\xi^s x$ for $s \in [a,b]$ with respect to the flow $\psi$. Among all points in $V_1$ equivalent to $y$ with respect to $\psi$, each is obtained from the corresponding $\xi^s x$ by $\psi^{-t}$ for a unique $t$ ($\psi^{-t}$ is non-singular with respect to $t$). Thus, $h$ has a one-to-one inverse. Moreover, $h$ is differentiable because it is the composition of differentiable functions. Thus, $h$ is a smooth bijection.

It follows then, by a change of variables and Fubini's theorem, that

$$(50) \qquad \mathcal{L}(V_1) = \int_{V_1} d\mathcal{L} = \int_a^b \int_I J(s,t) \, dt \, ds \leq Cr,$$

where $J$ is the absolute Jacobian determinant, and $C = \kappa |I| \sup_{s,t} J(s,t) < \infty$ by the compactness of $[a,b] \times I$.

Suppose now we take $\gamma$ to be a general curve as stipulated in the theorem, of length $r > 0$. Then, by taking $r$ sufficiently small, Lemma 1 shows that the gradients of $g_X$ along $\gamma$ point in the same direction up to $O(r)$. We can thus find a segment of an integral curve of $\xi$ within $O(r)$ of $\gamma$ that cuts across all the equivalence classes that $\gamma$ hits. The Lebesgue measure of $\mathcal{V}$ applied to the latter curve is bounded by $Cr$, and the additional area between the curves is $O(r^2)$.

Having proved the result for sufficiently small segments, a general curve can be decomposed into smaller segments for which the above arguments apply. The measure of the $\mathcal{V}$-induced set is then bounded above by a sum of the upper bounds along the curve, which gives a bound of $O(r)$ as was to be proved. □

**7. Proofs for Section 3.** Through this section, $\sup_x$ denotes the supremum over all regular points, hence, the rates of convergence resulting from Theorem 3 and Corollary 1 are uniform over the set of regular points of $g_X$.

To prove Theorem 3, we proceed as follows. Define

$$(51) \qquad p_n^*(x) = \frac{1}{n} \sum_{i=1}^n \frac{1}{\nu_n} K\left(\inf_{z \in P(X_i)} \frac{\|x - z\|}{\nu_n}\right).$$

Then,

$$|\widehat{p}_n(x) - p(x)| \leq |\widehat{p}_n(x) - p_n^*(x)| + |p_n^*(x) - \mathbb{E}[p_n^*(x)]| + |\mathbb{E}[p_n^*(x)] - p(x)|.$$

We bound these terms using the following three results.

THEOREM 7. *Under the assumptions in Theorem 3,*

$$\sup_x |\widehat{p}_n(x) - p_n^*(x)| = O\left(\sqrt{\frac{\log(1/h_n)}{nh_n^4}}\right) + O(h_n^2) \qquad a.s.$$

THEOREM 8. *Under the assumptions in Theorem 3,*

$$\sup_x |p_n^*(x) - \mathbb{E}[p_n^*(x)]| = O_P\left(\sqrt{\frac{\log n}{n\nu_n}}\right).$$

THEOREM 9. *Under the assumptions in Theorem 3,*

$$\sup_x |\mathbb{E}[p_n^*(x)] - p(x)| = O(\nu_n).$$

Theorem 3 follows by combining these results. Before proving these theorems, we need a few preliminary results.

LEMMA 3. *If (K1)–(K3) hold and $nh_n^2/\log(1/h_n) \to \infty$, then*

$$\sup_x |\widehat{g}_n(x) - g_X(x)| = O\left(\sqrt{\frac{\log h_n^{-1}}{nh_n^2}}\right) + O(h_n^2) \qquad a.s.$$

*If, in addition, $nh_n^4/\log(1/h_n) \to \infty$, then*

$$\sup_x |\nabla \widehat{g}_n(x) - \nabla g_X(x)| = O\left(\sqrt{\frac{\log h_n^{-1}}{nh_n^4}}\right) + O(h_n^2) \qquad a.s.$$

The first result above is Theorem 2.3 of Giné and Guillou (2002). The second result may be proved similarly; see also Giné and Koltchinskii (2006).

Let

$$D(x,y) = \inf_{s \in P(y)} \|s - x\|, \qquad \widehat{D}_n(x,y) = \inf_{s \in \widehat{P}(y)} \|s - x\|.$$

LEMMA 4.  *If* $nh_n^4 / \log(1/h_n) \to \infty$, *then*

$$\sup_x \sup_y |D(x,y) - \widehat{D}_n(x,y)| = O\left(\sqrt{\frac{\log h_n^{-1}}{nh_n^4}}\right) + O(h_n^2) \qquad a.s.$$

PROOF.  Follows from compactness and Lemma 3.  □

PROOF OF THEOREM 7.  Note that

$$\widehat{p}_n(x) - p_n^*(x) = \frac{1}{n\nu_n} \sum_{i=1}^n (K(\widehat{D}_n(x, X_i)/\nu_n) - K(D(x, X_i)/\nu_n)).$$

Since $K'$ is uniformly bounded, the result follows by Taylor expanding $K(\widehat{D}_n(x, X_i)/\nu_n)$ around $D(x, X_i)$ and applying the previous lemma.  □

Since $K$ is nonincreasing, there exists a cumulative distribution function $M$ with support $[0, \infty)$ such that

$$K(x) = \int_0^\infty \frac{1}{s} 1_{[0,s]}(x) \, dM(s).$$

Note that $K(x) = \int_x^\infty \frac{1}{s} dM(s)$ and that condition (K3) above implies that all moments of $M$ are finite and that $1 - M(x) = O(x^2 e^{-x})$.

The path density estimator $p_n^*$ in (51) can be written as

$$p_n^*(x) = \frac{1}{n} \sum_{i=1}^n \frac{1}{\nu_n} K\left(\frac{D(x, X_i)}{\nu_n}\right) = \frac{1}{n} \sum_{i=1}^n \frac{1}{\nu_n} \int_0^\infty \frac{1}{s} 1_{[0,s]}\left(\frac{D(x, X_i)}{\nu_n}\right) dM(s)$$

$$= \frac{1}{n\nu_n} \sum_{i=1}^n \int_0^\infty \frac{1}{s} 1_{[0,s\nu_n]}(D(x, X_i)) \, dM(s).$$

To prove Theorem 8 we will use Talagrand's (1994) inequality. The version we use is from Ginè and Guillou (2002). If $\Gamma$ is a class of functions, let $N(\Gamma, L_2(P), \delta)$ be the smallest number of balls of radius $\delta$ needed to cover $\Gamma$ with respect to the metric $\sqrt{\int (f(x) - g(x))^2 \, dP(x)}$, with $f, g \in \Gamma$.

THEOREM 10 (Ginè and Guillou). *Let $\Gamma$ be a class of uniformly bounded functions such that there exist $A \geq 3\sqrt{e}$ and $d \geq 1$ for which*

$$\sup_P N\left(\Gamma, L_2(P), \varepsilon\sqrt{\int G^2(x)\,dP(x)}\right) \leq \left(\frac{A}{\varepsilon}\right)^d, \tag{52}$$

*where $G$ is an envelope for $\Gamma$, and the supremum is over all probability measures. Let $\sigma^2 \geq \sup_{\gamma \in \Gamma} \mathrm{Var}(\gamma(X))$ and $U \geq \sup_{\gamma \in \Gamma} \|\gamma\|_\infty$ and $0 < \sigma \leq U$. Then, there exist constants $B, C$ and $L$ such that the following is true. If*

$$\varepsilon \geq \frac{C}{n}\left(U \log\left(\frac{AU}{\sigma}\right) + \sqrt{n}\sigma\sqrt{\log\left(\frac{AU}{\sigma}\right)}\right),$$

*then*

$$\mathbb{P}\left(\sup_{\gamma \in \Gamma}\left|\frac{1}{n}\sum_{i=1}^n \gamma(X_i) - \mathbb{E}(\gamma(X_i))\right| > \varepsilon\right)$$
$$\leq L\exp\left\{-\frac{n\varepsilon}{LU}\log\left(1 + \frac{n\varepsilon U}{L(\sqrt{n}\sigma + U\sqrt{\log(AU/\sigma)})^2}\right)\right\}.$$

PROOF OF THEOREM 8. Let $\gamma_{x,s} = 1_{[0,s\nu_n]}(D(x,\cdot))$ and write $p_n^*$ as

$$p_n^*(x) = \frac{1}{n\nu_n}\sum_{i=1}^n \int_0^\infty \frac{1}{s} 1_{[0,s\nu_n]}(D(x,X_i))\,dM(s)$$
$$= \int_0^\infty \frac{1}{s\nu_n}\frac{1}{n}\sum_{i=1}^n \gamma_{x,s}(X_i)\,dM(s).$$

Theorem 11 below shows that the class of functions

$$\Gamma_n = \{\gamma_{x,s}, x \in U_0, 0 \leq s \leq \log n\}$$

satisfies the covering condition (52), hence Talagrand's inequality holds with $U = 1$ and $\sigma^2 = O(\nu_n \log n)$. Set $Q_n = \sqrt{\log n}$. Then,

$$\mathbb{P}\left(\sup_{\gamma_{x,s} \in \Gamma_n}\left|\frac{1}{n}\sum_{i=1}^n \gamma_{x,s}(X_i) - \mathbb{E}[\gamma_{x,s}(X)]\right| > Q_n\sqrt{\frac{\nu_n \log n}{n}}\right)$$
$$\leq L\exp\left\{-\frac{nQ_n\sqrt{\nu_n \log n/n}}{L}\right.$$
$$\left.\times \log\left(1 + \frac{nQ_n\sqrt{\nu_n \log n/n}}{L(\sqrt{n\nu_n \log n} + \sqrt{\log(A/\sqrt{\nu_n \log n})})^2}\right)\right\} = o(1)$$

so that—with probability tending to 1—we have

$$\sup_{\gamma_{x,s} \in \Gamma_n}\left|\frac{1}{n}\sum_{i=1}^n \gamma_{x,s}(X_i) - \mathbb{E}[\gamma_{x,s}(X)]\right| = O\left(\sqrt{\frac{\nu_n \log n}{n}}\right).$$

It follows that

$$\sup_x |p_n^*(x) - \mathbb{E}[p_n^*(x)]|$$

$$\leq \sup_x \int_0^\infty \frac{1}{s\nu_n} \left| \frac{1}{n} \sum_{i=1}^n \gamma_{x,s}(X_i) - \mathbb{E}[\gamma_{x,s}(X)] \right| dM(s)$$

$$\leq \sup_x \int_0^{\log n} \frac{1}{s\nu_n} \left| \frac{1}{n} \sum_{i=1}^n \gamma_{x,s}(X_i) - \mathbb{E}[\gamma_{x,s}(X)] \right| dM(s) + \frac{1}{\nu_n} \int_{\log n}^\infty \frac{1}{s} dM(s)$$

$$\leq \sup_x \int_0^{\log n} \frac{1}{s\nu_n} \sup_{s \in (0, \log n)} \left| \frac{1}{n} \sum_{i=1}^n \gamma_{x,s}(X_i) - \mathbb{E}[\gamma_{x,s}(X)] \right| dM(s)$$

$$+ \frac{1}{\nu_n} K(\log n)$$

$$= \sup_x \sup_{s \in (0, \log n)} \left| \frac{1}{n} \sum_{i=1}^n \gamma_{x,s}(X_i) - \mathbb{E}[\gamma_{x,s}(X)] \right| \frac{1}{\nu_n} \int_0^{\log n} \frac{1}{s} dM(s)$$

$$+ \frac{1}{\nu_n} K(\log n)$$

$$\leq \frac{K(0)}{\nu_n} \sup_{\gamma_{x,s} \in \Gamma_n} \left| \frac{1}{n} \sum_{i=1}^n \gamma_{x,s}(X_i) - \mathbb{E}[\gamma_{x,s}(X)] \right| + \frac{K(\log n)}{\nu_n}.$$

Hence,

$$\sup_x |p_n^*(x) - \mathbb{E}[p_n^*(x)]| = \frac{K(0)}{\nu_n} O_P\left( \sqrt{\frac{\nu_n \log n}{n}} \right) + \frac{1}{\nu_n} O\left( \frac{\log n}{n} \right)$$

$$= O_P\left( \sqrt{\frac{\log n}{\nu_n n}} \right). \qquad \square$$

PROOF OF THEOREM 9. Note that $D(x, X_i) \leq s\nu_n$ if and only if $X_i \in \mathcal{V}(B(x, s\nu_n))$, hence

$$\mathbb{E}[1_{[0, s\nu_n]}(D(x, X))] = \mu_X(\mathcal{V}(B(x, s\nu_n))) = \pi(B(x, s\nu_n)).$$

The expected value of $p_n^*$ is

$$\mathbb{E}[p_n^*(x)] - p(x) = \frac{1}{\nu_n} \mathbb{E} \int_0^\infty \frac{1}{s} 1_{[0, s\nu_n]}(D(x, X)) \, dM(s) - p(x)$$

$$= \frac{1}{\nu_n} \int_0^\infty \frac{1}{s} \mathbb{E}[1_{[0, s\nu_n]}(D(x, X))] \, dM(s) - p(x)$$

$$= \frac{1}{\nu_n} \int_0^\infty \frac{1}{s} \pi(B(x, s\nu_n)) \, dM(s) - p(x)$$

$$= \frac{1}{\nu_n} \int_0^\infty \frac{1}{s}[\pi(B(x,s\nu_n)) - s\nu_n p(x)] \, dM(s)$$
$$= O(\nu_n).$$

To verify that the last term is $O(\nu_n)$, note that the integrand is not greater than $\frac{1}{s} + \nu_n p(x)$ and, from Theorem 2, there exists $r_x$ such that for $s\nu_n < r_x$ the integrand is smaller than $\frac{1}{s}Cs^2\nu_n^2$, so that

$$\left|\frac{1}{s}\pi(B(x,s\nu_n)) - s\nu_n p(x)\right| \leq \begin{cases} Cs\nu_n^2, & s < \frac{r_x}{\nu_n}, \\ \frac{\nu_n}{r_x} + \nu_n p(x), & s \geq \frac{r_x}{\nu_n}. \end{cases}$$

Also, $r \equiv \inf_x r_x > 0$. Hence

$$\left|\frac{1}{\nu_n} \int_0^\infty \frac{1}{s}[\pi(B(x,s\nu_n)) - s\nu_n p(x)] \, dM(s)\right|$$
$$\leq \int_0^{r/\nu_n} Cs\nu_n \, dM(s) + \int_{r/\nu_n}^\infty \left(\frac{1}{r} + p(x)\right) dM(s)$$
$$\leq C\nu_n \int_0^\infty s \, dM(s) + \left(\frac{1}{r} + p(x)\right)\left(1 - M\left(\frac{r}{\nu_n}\right)\right),$$

which is $O(\nu_n)$, due to the tail condition on $M$ and the existence of its moments. Thus,

$$\sup_x |\mathbb{E}[p_n^*(x)] - p(x)| = O(\nu). \qquad \square$$

THEOREM 11. *Let $k$ be the number of critical points of $g_X$ and let $m \leq k$ be the number of modes. The VC dimension of the class of sets*

$$\mathcal{R} = \{\mathcal{V}(B(x, s\nu_n)) : x \in U_0, 0 \leq s \leq \log n\}$$

*is less than $4m$. Hence, the covering condition (52) holds.*

PROOF. Let $b_1, \ldots, b_m$ denote the modes of $g_X$ and let $\mathcal{E}(x) = \{\psi^t(x), 0 \leq t \leq \infty\}$. Note that $\mathcal{E}(x)$ is a smooth curve starting at $x$ and ending at some mode of $g_X$. Partition $U_0$ into sets $U_1, \ldots, U_m$ where $x \in U_j$ if $\psi^\infty(x) = b_j$. Let $F$ be a finite set containing $4m$ points. Then there exists at least one $U_j$ such that $F \cap U_j$ has at least 4 points. Let $G \subset F \cap U_j$ be a subset of size 4. So $G = \{x_1, x_2, x_3, x_4\}$, say. We will show that $G$ cannot be shattered.

*Step* 1. If $x_k \in \mathcal{E}(x_\ell) - \{b_j\}$ for $x_k \neq x_\ell$ (where $x_k, x_\ell \in G$) then clearly $G$ cannot be shattered. Thus we can assume that $\mathcal{E}_k \cap \mathcal{E}_\ell = \{b_j\}$ for each $k \neq \ell$.

*Step* 2. Let $R = \mathcal{V}(B(z, \rho)) \in \mathcal{R}$. We claim that $R$ picks out $A \subset G$ if and only if

(53) $\qquad\qquad B(z, \rho) \cap \mathcal{E}(x_j) \neq \varnothing \qquad \text{for each } x_j \in A$

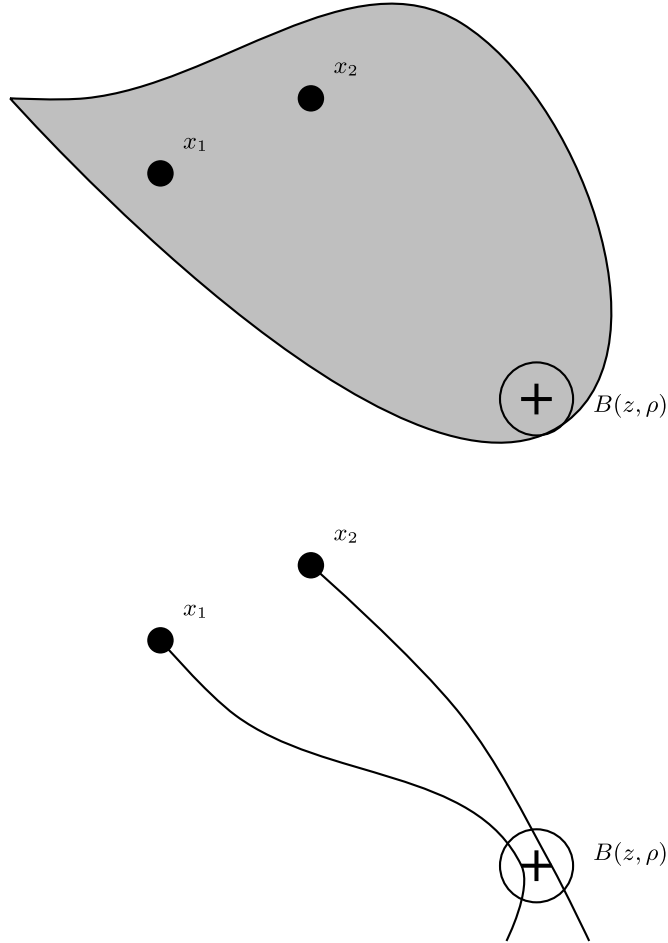

FIG. 5. *The set $R = \mathcal{V}(B(z,\rho))$ picks out $\{x_1, x_2\}$ (above) if and only if the forward evolutions $\mathcal{E}(x_1)$ and $\mathcal{E}(x_2)$ hit the ball $B(z,\rho)$ (below).*

and

(54) $$B(z,\rho) \cap \mathcal{E}(x_j) = \varnothing \qquad \text{for each } x_j \in G - A.$$

This follows since $x \in \mathcal{V}(B(z,\rho))$ if and only if $B(x,\rho) \cap \mathcal{E}(x) \neq \varnothing$. See Figure 5.

*Step* 3. Place a small ball around $b_j$ and renumber the points in $G$ by the angle of the vector from $z$ to the intersection of $\mathcal{E}(x_j)$ with the boundary of the ball. See Figure 6.

*Step* 4. We claim that if there exists an $R$ that picks out $A = \{x_1, x_3\}$ then there does not exist an $R'$ that picks out $\{x_2, x_4\}$. Suppose $R = B(z,\rho)$ picks out $A$. Because of (53) and (54), both $\mathcal{E}(x_1)$ and $\mathcal{E}(x_3)$ pass through $B(z, \rho)$ but neither $\mathcal{E}(x_2)$ and $\mathcal{E}(x_4)$ pass through $B(z, \rho)$. See Figure 7. Choose $y_1 \in \mathcal{E}(x_1) \cap B(z, \rho)$ and $y_2 \in \mathcal{E}(x_3) \cap B(z, \rho)$. Define a closed curve $C$ as follows. $C$ begins at the mode $b_j$, follows $\mathcal{E}(x_1)$ to $y_1$, connects $y_1$ to $y_2$ by any smooth curve contained in $B(z, \rho)$ and the follows $\mathcal{E}(x_2)$ back to $b_j$. Now $C$ encloses $\mathcal{E}(x_2)$ and excludes $\mathcal{E}(x_4)$. Hence, there is no ball $B(z', \rho)$ satisfying (53) and (54) for $A' = \{x_2, x_4\}$. Thus, there does not exist an $R'$ that picks out $\{x_2, x_4\}$. We conclude that $G$ cannot be shattered. □

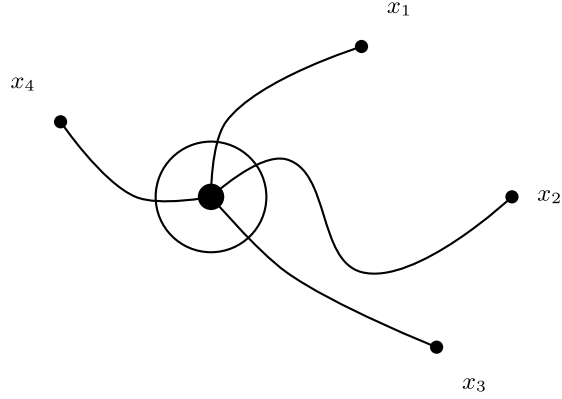

FIG. 6. *The curves $\mathcal{E}(x_1), \ldots, \mathcal{E}(x_4)$ can be ordered by their intersection with a small disc around $b_j$.*

### 8. Proofs for Section 4.

PROOF OF THEOREM 4. We prove the theorem in the background-free case. The result follows in the general case because the gradient of $g_X$ on $U_0$ changes by a constant factor, yielding the same function $p(x)$. Let $D(x, r)$ denote the set $\mathcal{V}(\gamma([\alpha_r, \beta_r]))$ in Theorem 6(4b) within the ball of radius $r$ around $x$.

Consider the following intermediate claims:

1. $x \preceq y$ implies $g_X(x) \leq g_X(y)$.

Let $\gamma(t) = \psi^{-t} y$, $t \geq 0$. Then, $\gamma'(t) = -\nabla g_X(\gamma(t))$ and there exists $t_1 > 0$ such that $\gamma(t_1) = x$. Let $v(t) = g_X(\gamma(t))$. Then, $v(0) = g_X(y)$, $\dot{v}(t) = \nabla g_X(\gamma(t)) \cdot \gamma'(t) = -\|\nabla g_X(\gamma(t))\|^2 < 0$. Hence, $v(t) = v(0) + \int_0^t \dot{v}(t)\, dt \leq v(0)$, which proves the claim.

2. $g_X(x) \leq \varphi_\sigma(d(x, \mathcal{A}))$.

Let $\widetilde{x}$ denote the point in $\mathcal{A}$ that is closest to $x$. Then, $g_X(x)$ is no greater than $g(x)$, where $g$ is the normal mixture with all its mass at $\widetilde{x}$. Thus, $g(x) = \varphi_\sigma(x - \widetilde{x}) = \varphi_\sigma(d(x, \mathcal{A}))$, proving the claim.

3. If $x \in B(\mathcal{A}, d(\lambda) + \varepsilon)^c$ and $r < \varepsilon$, every $y \in \mathcal{V}(B(x, r))$ and every $y \in \mathcal{V}(D(x, r))$ satisfies

(55) $$g_X(y) \leq \varphi_\sigma(d(\lambda)).$$

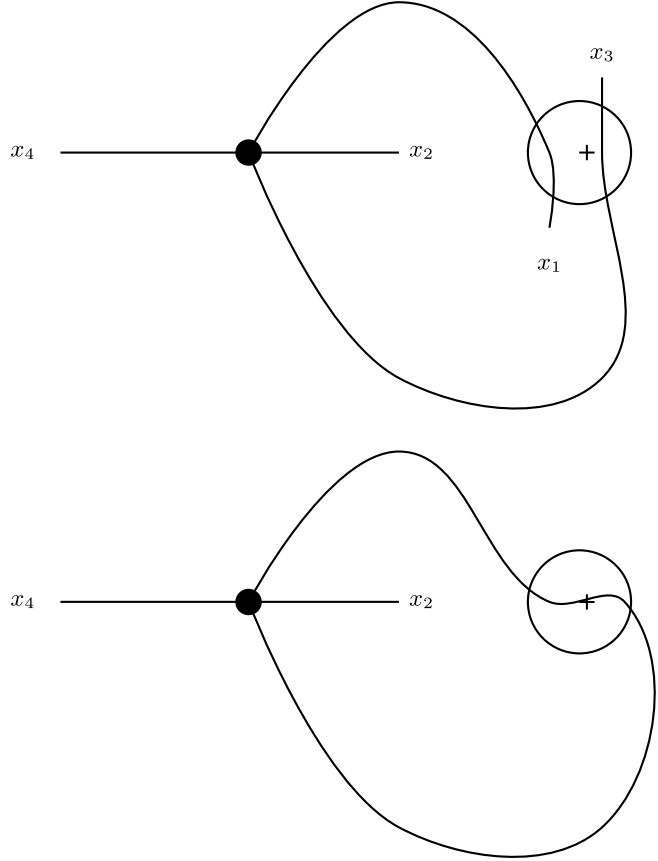

Fig. 7.  *If $\mathcal{V}(B(z,\rho))$ picks out $\{x_1, x_3\}$ then $\mathcal{E}(x_1)$ and $\mathcal{E}(x_2)$ hit the ball $B(z,\rho)$ (top). We can then join the two curves $\mathcal{E}(x_1)$ and $\mathcal{E}(x_2)$ forming a closed curve $C$ that isolates $x_2$ from $x_4$. Note that all $\mathcal{E}(x_j)$ must end at the mode (the large dot) since these points are all in the same element of the partition $U_j$.*

Every element of $B(x,r)$ is farther than $d(\lambda)$ from $\mathcal{A}$, where $d(\lambda)$ is defined in (14), so the first claim follows by 1 and 2. By the construction in the proof of Theorem 6, $D(x,r) \subset B(x,r)$, so the second claim follows as well.

4. For regular points $x$ of $g_X$, $\mathcal{L}(\mathcal{V}(D(x,r))) \leq Cr + O(r^2)$, for a constant that depends only on $U_0$.

The claim follows pointwise directly from Lemma 2. The proof of that lemma shows that because $U_0$ is compact, the constants for each point are bounded.

Now we prove the main result. If $x \in \mathcal{M}$, the local maxima of $g_X$, then $p(x) = \infty$ by Theorem 6. If $x \in \mathcal{H}$, then $p(x) \leq 4\sup_y p(y)$, where the supremum is over any sufficiently small ball around $x$. Hence, if $x \in B(\mathcal{A}, d(4\lambda) +$ $\varepsilon)^c$, the regular points in a small enough ball around $x$ will lie in $B(\mathcal{A}, d(4\lambda) + \varepsilon)^c$.

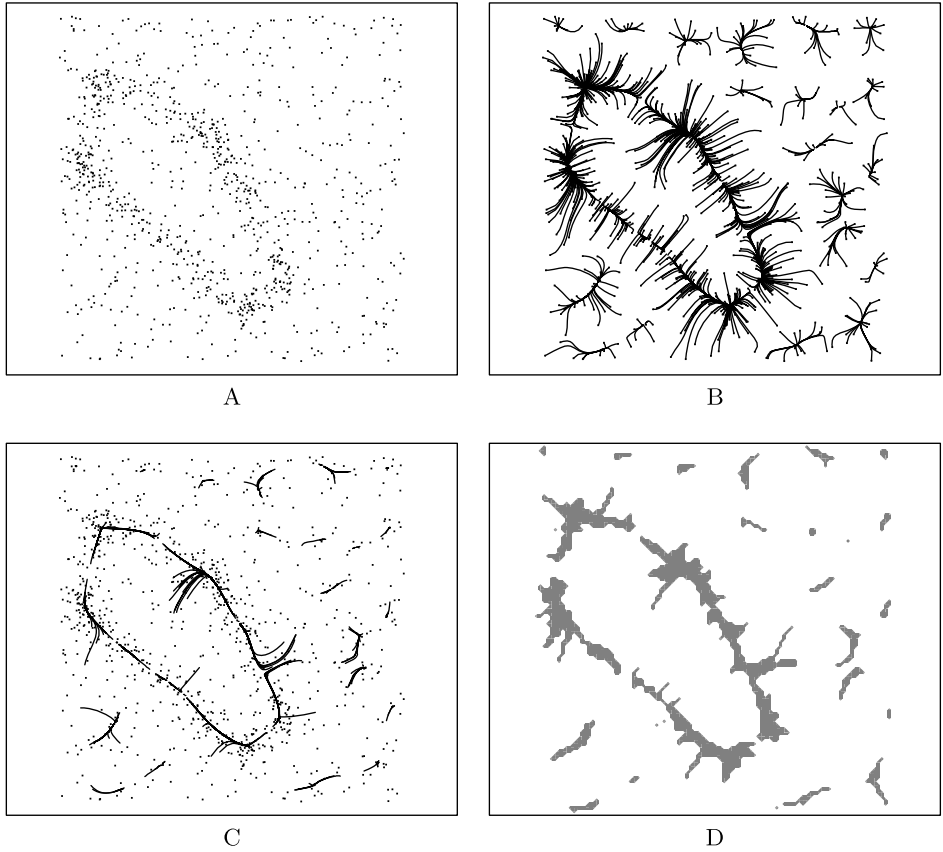

FIG. 8. *Simulated example with background noise: data set (plot* A*), steepest ascent paths of all data points (plot* B*), paths after cutting the first seven iterations (plot* C*) and level set at 90% quantile of the estimated path density (plot* D*).*

If $x \in B(\mathcal{A}, d(\lambda) + \varepsilon)^c$ is a regular point of $g_X$, then for small $r > 0$, $\mathcal{L}(\mathcal{V}(D(x,r))) \leq Cr + O(r^2)$, for a constant $C$ independent of $x$. Hence,

(56) $$\mu_X(\mathcal{V}(D(x,r))) \leq Cr\phi_\sigma(d(\lambda)) + O(r^2).$$

Thus $p(x) \leq C\phi_\sigma(d(\lambda)) = \lambda$, and the result follows. $\square$

PROOF OF THEOREM 5. The theorem follows directly from the previous results and from Theorem 3 of Cuevas and Fraiman (1997). $\square$

C. R. Genovese
L. Wasserman
Department of Statistics
Carnegie Mellon University
Pittsburgh, Pennsylvania 15213
USA
E-mail: genovese@stat.cmu.edu
       larry@stat.cmu.edu

M. Perone-Pacifico
Dipartimento di Statistica
Sapienza Università di Roma
Piazzale A. Moro 5
00183 Roma
Italy
E-mail: marco.peronepacifico@uniroma1.it

I. Verdinelli
Dipartimento di Statistica
Sapienza Università di Roma
Piazzale A. Moro 5
00183 Roma
Italy
and
Carnegie Mellon University
5000 Forbes Avenue
Pittsburgh, Pennsylvania 15213
USA
E-mail: isabella@stat.cmu.edu